\RequirePackage{fix-cm}

\documentclass[smallextended]{svjour3}

\makeatletter
\let\cl@chapter\@empty
\makeatother

\smartqed

\usepackage{graphicx}
\usepackage{amssymb}
\usepackage{amsmath}
\usepackage{xcolor}
\usepackage{tikz}
\usepackage{caption}
\usepackage{subcaption}
\usepackage{algorithm}
\usepackage{algorithmic}
\usepackage{pgfplots}
\usepackage{url}
\usepackage{comment}

\usepackage{cleveref}

\usetikzlibrary{patterns}
\def\amen{{{\sc AMEn}}}

\newcommand*{\bigtimes}{\mathop{\raisebox{-.5ex}{\hbox{\huge{$\times$}}}}} 

\definecolor{mycolor1}{rgb}{0.00000,0.44700,0.74100}%
\definecolor{mycolor2}{rgb}{0.85000,0.32500,0.09800}%
\definecolor{mycolor3}{rgb}{0.92900,0.69400,0.12500}%
\definecolor{mycolor4}{rgb}{0.49400,0.18400,0.55600}%
\definecolor{mycolor5}{rgb}{0.46600,0.67400,0.18800}%

\begin{document}

\title{Projection-based low-rank assembly in IgA}
\subtitle{}


\author{Tom-Christian Riemer         \and
        Martin Stoll 
}


\institute{T.-C. Riemer \at
                Technical University Chemnitz, Department of Mathematics, Reichenhainer Str. 39, 09126 Chemnitz, Saxony, Germany \\
                Tel.: +49 371 531 34059 \\
                \email{tom-christian.riemer@math.tu-chemnitz.de}           
                \and
                M. Stoll \at
                Technical University Chemnitz, Department of Mathematics, Reichenhainer Str. 39, 09126 Chemnitz, Saxony, Germany \\
                Tel.: +49 371 531 33705 \\
                \email{martin.stoll@math.tu-chemnitz.de}          
}

\date{Received: date / Accepted: date}

\maketitle

\begin{abstract}
Isogeometric Analysis (IgA) uses the same spline functions for the representation of the computational domain and for the approximation of the solution, which allows exact geometry descriptions but leads to mass and stiffness matrices that are expensive to assemble and to store, especially in three dimensions. In this paper, we present a projection-based low-rank approach for assembling the mass tensor and the stiffness tensor of orientation-preserving tensor-product B-spline geometries. For the mass tensor, we exploit the polynomial structure of the determinant of the Jacobian of the geometry map and represent it in reduced spline product spaces by univariate coefficient transfer operators. Using exact quadrature and without tensor truncation, the resulting low-rank tensor is an exact reformulation of the standard Galerkin mass tensor. For the stiffness tensor, we separate the rational weight function into a polynomial numerator, again represented in reduced spline product spaces, and the reciprocal determinant, which is not a spline function in general and is therefore approximated by an $L^2$-projection onto a tensor-product spline space. Both constructions are carried out entirely in the tensor-train (TT) format, where the projection system is solved by the alternating minimal energy (\amen) method, so that full high-order coefficient tensors are never formed and the multidimensional integrals reduce to univariate integrals and contracted products. The method is implemented in \textsc{MATLAB} using the \textsc{GeoPDEs} and \textsc{TT-Toolbox} packages. Our numerical experiments show that the proposed method is competitive with the full assembly and with the interpolation-based low-rank method. It is applicable in two situations in which interpolation is problematic, namely when the interpolation system becomes singular and for nearly singular geometries. The construction is restricted to orientation-preserving tensor-product B-spline geometries and does not cover NURBS parameterizations.
\end{abstract}

\section{Introduction}
\label{sec:Introduction}

Isogeometric Analysis (IgA), introduced in~\cite{CAD}, is a discretization technique for the numerical solution of partial differential equations that aims to bridge the gap between computer-aided design and finite element analysis. Its central idea is to use the same spline functions, B-splines or non-uniform rational B-splines (NURBS), for the representation of the computational domain and for the approximation of the solution, so that exact or highly accurate geometry descriptions can be used directly in the analysis. The tensor-product structure and the high degrees of the basis functions, however, make the resulting mass and stiffness matrices expensive to assemble and to store. This is especially critical in three dimensions, where multivariate quadrature is costly and the overlap of basis functions leads to large and dense matrix structures. Efficient assembly techniques and compact matrix representations are therefore essential to make IgA competitive for large-scale problems.

Low-rank tensor methods have become a promising tool in this context. The main idea is to exploit the tensor-product structure of multivariate B-spline spaces and to represent the geometry-dependent weight functions in separated form, so that the arising multidimensional integrals reduce to sums of products of univariate integrals. This idea was introduced in~\cite{angelos1} and later combined with tensor-train (TT) techniques in~\cite{BuengerDolgovStoll:2020}, where the coefficient tensors are represented in TT format and computed without forming the full tensors explicitly. Low-rank techniques have also been used directly in the solution process of isogeometric linear systems, cf.~\cite{montardini2023lowsp}, and have since been extended to multi-patch geometries, cf.~\cite{montardini2024lowrank,Riemer2025}, to (truncated) hierarchical B-splines, cf.~\cite{Riemer2026}, and to matrix-free inexact preconditioning, cf.~\cite{Mika2026}. These developments show that low-rank methods remain an active research direction, both for matrix assembly and for the efficient solution of the resulting linear systems.

In this paper, we present a new low-rank approach for assembling the mass tensor and the stiffness tensor for orientation-preserving tensor-product B-spline geometries. In contrast to interpolation-based approaches, the mass tensor assembly exploits the polynomial structure of the determinant of the Jacobian of the geometry map, which is represented in a suitable reduced spline product space by univariate coefficient transfer operators. If the required coefficients and univariate integrals are computed with exact quadrature and no tensor truncation is applied, the resulting low-rank mass tensor is an exact reformulation of the standard Galerkin mass tensor. For the stiffness tensor, the matrix-valued weight function contains the inverse of the Jacobian of the geometry map and is rational in general. We separate it into a polynomial numerator part, which is again represented in reduced spline product spaces exactly up to quadrature and truncation, and the reciprocal determinant, which is in general not a spline function and is therefore approximated by a projection onto a suitable tensor-product spline space. Both constructions are carried out in the TT format, where the projection system for the reciprocal determinant is solved by the alternating minimal energy (\amen) method~\cite{amen}, so that full high-order coefficient tensors are never formed. The multidimensional integrals reduce then to univariate integrals and contracted products. We point out that the construction is restricted to orientation-preserving tensor-product B-spline geometries and does not cover NURBS parameterizations.

Outside of providing a separated representation, the proposed constructions are of particular advantage in two situations in which the interpolation-based assembly is problematic. First, the proposed mass and stiffness assembly remain applicable when the interpolation system becomes singular, which can occur and then prevents the interpolation-based assembly from producing a usable low-rank representation. Second, for geometry maps that are almost singular, that is, whose determinant of its Jacobian nearly vanishes in part of the parametric domain, the projection-based stiffness assembly shows a slight but consistent accuracy advantage over the interpolation variant. Both effects are confirmed by the numerical experiments in \cref{sec:Numerics}.

All numerical experiments are implemented in \textsc{MATLAB}. We use the \textsc{GeoPDEs} package~\cite{geopdes3.0,geopdes_new} to set up the isogeometric discretization and to assemble reference mass and stiffness matrices in full format, the \textsc{TT-Toolbox}~\cite{tt-toolbox} for representing tensors in TT format and performing the corresponding operations, and, for comparison, an adapted version of the interpolation-based implementation from~\cite{BuengerCode}.

The paper is organized as follows. \cref{sec:Preliminaries} recalls orientation-preserving tensor-product B-spline geometries and the interpolation-based low-rank assembly. \cref{sec:Projection} develops the projection-based construction of the low-rank mass and stiffness tensors. \cref{sec:Limitations} discusses the limitations of the approach and the implementation details, and \cref{sec:Numerics} reports the numerical experiments before we conclude.

\section{Preliminaries}
\label{sec:Preliminaries}

\subsection{Orientation preserving tensor-product B-spline geometries}
\label{sec:Tensor-ProductGeometries}
Let $p\in\mathbb N_0$ and let $\Xi= \left[ \xi_1, \ldots, \xi_{n+p+1} \right]$ be a knot vector with $0=\xi_1\le \cdots \le \xi_{n+p+1}=1$. The univariate B-spline basis $\mathcal{B}_p  \left(  \Xi  \right)  = \left\{ \beta_{i,p} \colon [0,1] \to [0,1] \right\}_{i=1}^n$ can be defined by the Cox-de Boor recursion, cf.~\cite{cox:b_splines,DEBOOR197250}, where quotients with zero denominator are set to zero. We denote by $\mathbb{S}_p  \left(  \Xi  \right)  = \operatorname{span}  \left(  \mathcal{B}_p  \left(  \Xi  \right)   \right)$ the linear space spanned by the univariate basis $\mathcal{B}_p  \left(  \Xi  \right)$ and omit the degree $p$ and the dependence on $\Xi$ whenever they are clear from the context. In this paper, we consider B-splines associated with open knot vectors, in which the first and the last knot value are repeated $p+1$ times, i.e. $\xi_1 = \ldots = \xi_{p+1} = 0$ and $\xi_{n+1} = \ldots = \xi_{n+p+1} = 1$. Each B-spline has support $\operatorname{supp}\left(\beta_{i,p}\right)=\left[\xi_i,\xi_{i+p+1}\right]$ and overlaps with at most $2p+1$ basis functions, including itself.

The derivative of a B-spline of degree $p$ is a linear combination of B-splines of degree $p-1$
\begin{equation}
        \label{eq:BSplineDerivative}
        \frac{\mathrm d}{\mathrm d\hat{x}} \beta_{i,p} \left(  \hat{x}  \right)  = \frac{p}{\xi_{i+p}-\xi_i} \beta_{i,p-1} \left(  \hat{x}  \right)  - \frac{p}{\xi_{i+p+1}-\xi_{i+1}} \beta_{i+1,p-1} \left(  \hat{x}  \right) ,
\end{equation}
again with zero quotients interpreted as zero. If $\Xi= \left[ \xi_1, \ldots, \xi_{n+p+1} \right]$ is the knot vector of the basis of degree $p$, then the corresponding reduced knot vector for the derivative space is $\Xi^\prime = \left[ \xi_2, \ldots, \xi_{n+p} \right]$. This means that for a spline function $u_h = \sum_{i=1}^{n} u_i \beta_{i,p}$, its derivative can be represented as
\begin{equation*}
        u^\prime_h = \sum_{i=1}^{n-1} p \, \frac{\left( u_{i+1}-u_i \right)}{\xi_{i+p+1}-\xi_{i+1}} \, \beta_{i,p-1},
\end{equation*}
where $\{\beta_{i,p-1}\}_{i=1}^{n-1}$ denotes the B-spline basis of degree $p-1$ associated with the reduced knot vector $\Xi^\prime$.

Let now $D \in\mathbb N$ denote the parametric dimension. For each $d=1,\ldots,D$, let $\Xi^{\left( d \right)}$ be a knot vector, let $p^{\left( d \right)} \in \mathbb N_0$ be the corresponding degree. Let $\mathcal{B}_{p^{\left( d \right)}} \left( \Xi^{\left( d \right)} \right)  = \left\{ \beta^{\left( d \right)}_{i^{\left( d \right)},p^{\left( d \right)}} \right\}_{i^{\left( d \right)}=1}^{n^{\left( d \right)}}$ be the associated univariate B-spline basis and the corresponding univariate spline space is then denoted by $ \mathbb{S}_{p^{\left( d \right)}} \left( \Xi^{\left( d \right)} \right)  = \operatorname{span}  \left(  \mathcal{B}_{p^{\left( d \right)}} \left( \Xi^{\left( d \right)} \right)   \right)$. Multivariate B-splines on the parameter domain $[0,1]^D$ are defined by tensor products of the univariate B-splines. For that we introduce a multi-index $\boldsymbol{i} =  \left( i^{\left( 1 \right)},\ldots,i^{\left( D \right)} \right)$ and $p =  \left( p^{\left( 1 \right)},\ldots,p^{\left( D \right)} \right)$, we set
\begin{equation*}
        \beta_{\boldsymbol{i},p}  \left(  \hat{x}  \right)  = \prod_{d=1}^D \beta^{\left( d \right)}_{i^{\left( d \right)},p^{\left( d \right)}}  \left(  \hat{x}^{\left( d \right)}  \right) , \quad \hat{x} =  \left( \hat{x}^{\left( 1 \right)},\ldots,\hat{x}^{\left( D \right)} \right) \in [0,1]^D.
\end{equation*}
Equivalently, if the univariate basis functions in $d$ are collected in the vector
\begin{equation*}
        b^{\left( d \right)}  \left(  \hat{x}^{\left( d \right)}  \right)  = \left[ \beta^{\left( d \right)}_{1,p^{\left( d \right)}} \left(  \hat{x}^{\left( d \right)}  \right) , \ldots, \beta^{\left( d \right)}_{n^{\left( d \right)},p^{\left( d \right)}} \left(  \hat{x}^{\left( d \right)}  \right)  \right]^\top \in \mathbb{R}^{n^{\left( d \right)}},
\end{equation*}
then all multivariate B-splines can be collected in the tensor
\begin{equation*}
        \mathbf{B}  \left(  \hat{x}  \right)  = \bigotimes^{D}_{d = 1} b^{\left( d \right)} \left(  \hat{x}^{\left( d \right)}  \right)  \in \mathbb{R}^{\left( n^{\left( 1 \right)},\ldots,n^{\left( D \right)} \right)} = \mathbb{R}^{\mathbf{n}},
\end{equation*}
with $\mathbf{n} =  \left( n^{\left( 1 \right)},\ldots,n^{\left( D \right)} \right)$. Its entries are the multivariate basis functions, 
\begin{equation*}
        \left(  \mathbf{B}  \left(  \hat{x}  \right)   \right) _{\left( i^{\left( 1 \right)},\ldots,i^{\left( D \right)} \right)} = \beta_{\boldsymbol{i},p}  \left(  \hat{x}  \right) .
\end{equation*}
The corresponding multivariate spline space is the tensor-product space
{\small
\begin{equation}
        \label{eq:tensor_product_space}
        \begin{aligned}
                \mathbb{S}_{p}  \left(  \Xi^{\left( 1 \right)},\ldots,\Xi^{\left( D \right)}  \right)  & = \operatorname{span}  \left(  \left\{ \beta_{\boldsymbol{i},p} \colon \; \boldsymbol{i} = \left( i^{\left( 1 \right)},\ldots,i^{\left( D \right)} \right)  \in \mathbf{I} = \bigtimes^{D}_{d = 1} \left\{1, \ldots, n^{\left( d \right)} \right\} \right\} \right) \\
                &  = \bigotimes_{d=1}^D \mathbb{S}_{p^{\left( d \right)}}  \left(  \Xi^{\left( d \right)}  \right) .
        \end{aligned}
\end{equation}
}
In particular, a discrete function $u_h\in \mathbb{S}_{p}  \left( \Xi^{\left( 1 \right)},\ldots,\Xi^{\left( D \right)} \right)$ can be written as
\begin{equation}
        \label{eq:DiscreteFunction}
        u_h \left(  \hat{x}  \right) = \sum_{\mathbf{i} \in \mathbf{I}} u_{\mathbf{i}} \beta_{\mathbf{i}}  \left(  \hat{x}  \right)   =  \left\langle \mathbf{u}, \mathbf{B}  \left(  \hat{x}  \right)   \right\rangle_F
\end{equation}
where the coefficients form a tensor $\mathbf{u} \in \mathbb{R}^{\mathbf{n}}$ and $ \left\langle \cdot, \cdot  \right\rangle_F$ represents the Frobenius product. 

There are several refinement strategies that can be applied to enrich a B-spline basis. There is $h$-refinement, which is realized by knot insertion. Here, additional knots are inserted into the knot vector, thereby increasing the number $n^{(d)}$ of basis functions without changing the spline degree. And there is $p$-refinement, where the degree $p^{(d)}$ of the B-splines is increased. In order to preserve the continuity of the original spline space, the multiplicities of the existing knots have to be increased accordingly~\cite{piegl,HugesIGA}. 

In IgA, spline functions are used to represent the computational domain $\Omega\subset\mathbb{R}^3$. Let $\mathcal{B} = \left\{ \beta_{\mathbf{i}} \colon [0,1]^3 \to \mathbb{R} \colon \; \mathbf{i}\in\mathbf{I} \right\}$ be a basis of a three-variate spline space $\mathbb{S}=\operatorname{span} \left( \mathcal{B} \right)$ (not necessarily a tensor-product basis and space). The geometry map is then defined by
\begin{equation*}
        G \colon [0,1]^3 \to \Omega, \qquad G \left( \hat{x} \right)  = \sum_{\mathbf{i}\in\mathbf{I}} \mathbf{C}_{\mathbf{i}} \beta_{\mathbf{i}} \left( \hat{x} \right) ,
\end{equation*}
where $\mathbf{C}_{\mathbf{i}}\in\mathbb{R}^3$, $\mathbf{i}\in\mathbf{I}$, are the so-called control points. We assume that the control points are chosen such that $G$ is at least injective. In particular, since we are interested in orientation-preserving geometries, we assume
\begin{equation}
        \label{eq:OrientationPreserving}
        \det  \left(  \nabla G \left( \hat{x} \right)   \right) > 0, \quad \hat{x}\in\ [0,1]^3.
\end{equation}

Very general spline spaces can be used to represent the computational domain $\Omega$, and these spaces do not necessarily possess a global tensor-product structure. For instance, local refinement techniques can be employed to accurately represent geometric details, see e.g.~\cite{GIANNELLI2016337}. In this paper, however, it is essential that the spline space has tensor-product structure, i.e. $\mathbb{S} = \mathbb{S}_{p}  \left( \Xi^{\left( 1 \right)},\Xi^{\left( 2 \right)},\Xi^{\left( 3 \right)}  \right)$ as in \cref{eq:tensor_product_space}. In the tensor-product case, the control points can then be collected in a tensor $\mathbf{C} \in \mathbb{R}^{{\left(3, n^{\left( 1 \right)}, n^{\left( 2 \right)}, n^{\left( 3 \right)} \right)}}$. More precisely, for $\mathbf{i}= \left( i^{\left( 1 \right)},i^{\left( 2 \right)},i^{\left( 3 \right)} \right) \in\mathbf{I}$, we write $\mathbf{C}_{\left( :,i^{\left( 1 \right)},i^{\left( 2 \right)},i^{\left( 3 \right)} \right)} = \mathbf{C}_{\mathbf{i}} \in \mathbb{R}^3$. As a result, the first mode of $\mathbf{C}$ corresponds to the physical coordinate dimension, while the remaining three modes correspond to the tensor-product B-spline basis. With this notation, the geometry map can be written as the contracted product
\begin{equation}
        \label{eq:GeometryMapTensor}
        G \left( \hat{x} \right)  = \mathbf{C} \cdot \mathbf{B} \left( \hat{x} \right)  =
        \begin{bmatrix}
                \left\langle \mathbf{C}^{(1)}, \ \mathbf{B} \left( \hat{x} \right)   \right\rangle_F \\
                \left\langle \mathbf{C}^{(2)}, \ \mathbf{B} \left( \hat{x} \right)   \right\rangle_F \\
                \left\langle \mathbf{C}^{(3)}, \ \mathbf{B} \left( \hat{x} \right)   \right\rangle_F
        \end{bmatrix}
        \in \mathbb{R}^3,
\end{equation}
where we refer to $\mathbf{C}^{(a)} = \mathbf{C}_{\left( a,:,:,: \right)} \in \mathbb{R}^{\left( n^{\left( 1 \right)}, n^{\left( 2 \right)}, n^{\left( 3 \right)} \right)}$, $a=1,2,3$, as coordinate control-point tensors from now on. Here, $\cdot$ denotes the contracted product over the last three modes of $\mathbf{C}$ and all modes of $\mathbf{B} \left( \hat{x} \right)  \in \mathbb{R}^{\left( n^{\left( 1 \right)}, n^{\left( 2 \right)}, n^{\left( 3 \right)} \right)}$. Equivalently, for $a=1,2,3$ the $a$-th component of the geometry map is given by
\begin{equation*}
        G_a \left( \hat{x} \right)  =  \left\langle \mathbf{C}^{(a)}, \mathbf{B} \left( \hat{x} \right)   \right\rangle_F = \sum_{i^{\left( 1 \right)}=1}^{n^{\left( 1 \right)}} \sum_{i^{\left( 2 \right)}=1}^{n^{\left( 2 \right)}} \sum_{i^{\left( 3 \right)}=1}^{n^{\left( 3 \right)}} \mathbf{C}_{\left( a,i^{\left( 1 \right)},i^{\left( 2 \right)},i^{\left( 3 \right)} \right)} \beta_{\left( i^{\left( 1 \right)},i^{\left( 2 \right)},i^{\left( 3 \right)} \right)}  \left( \hat{x} \right).
\end{equation*}

In addition to B-spline geometries, computational domains are often represented using NURBS which constitute a rational extension of B-splines and allow the exact representation of many geometries that cannot be represented exactly by polynomial B-splines, including conic sections. Because of their rational structure, multivariate NURBS basis functions are, in general, not separable into tensor products of univariate functions. As a consequence, the direct tensor-product structure exploited by the proposed method is lost and in the following, we restrict ourselves to B-spline geometries whose geometry maps are constructed from tensor-product B-spline bases.

\subsection{Low-rank methods in {IgA} based on interpolation}
\label{subsection:LowRankInterpolation}
Let $\hat{\mathbb{S}}_{\hat{p}}  \left(  \hat{\Xi}^{\left( 1 \right)}, \hat{\Xi}^{\left( 2 \right)}, \hat{\Xi}^{\left( 3 \right)}  \right)$ be a tensor-product solution space obtained from the geometry space $\mathbb{S}_{p} \left( \Xi^{\left( 1 \right)},\Xi^{\left( 2 \right)},\Xi^{\left( 3 \right)} \right)$ by refinement, i.e.
\begin{equation*}
        \mathbb{S}_{p} \left( \Xi^{\left( 1 \right)},\Xi^{\left( 2 \right)},\Xi^{\left( 3 \right)} \right)  \subseteq \hat{\mathbb{S}}_{\hat{p}}  \left(  \hat{\Xi}^{\left( 1 \right)}, \hat{\Xi}^{\left( 2 \right)}, \hat{\Xi}^{\left( 3 \right)}  \right),
\end{equation*}
and let $\left\{\hat{\beta}_{\hat{\mathbf{i}}} \right\}_{\hat{\mathbf{i}}\in\hat{\mathbf{I}}}$ denote its tensor-product basis with $\hat{\mathbf{I}} = \bigtimes_{d=1}^{3} \left\{1,\ldots,\hat{n}^{\left( d \right)}\right\}$. The corresponding physical solution space is
\begin{align*}
        V_h & = \left\{ \hat{v}_h \circ G^{-1} \colon \Omega \to \mathbb{R} \mid \hat{v}_h \in \hat{\mathbb{S}}_{\hat{p}}  \left(  \hat{\Xi}^{\left( 1 \right)}, \hat{\Xi}^{\left( 2 \right)},\hat{\Xi}^{\left( 3 \right)}  \right)  \right\} \\
        & = \operatorname{span} \left(  \left\{ \hat{\beta}_{\hat{\mathbf{i}}} \circ G^{-1} \colon \Omega \to \mathbb{R} \mid \hat{\mathbf{i}} \in \hat{\mathbf{I}} \right\}  \right).
\end{align*}

Following the Galerkin method and applying the coordinate transformation from $\Omega$ to $[0,1]^3$, the entries of the mass and stiffness matrices or tensors are given by
\begin{equation}
        \label{eq:MassAndStiffness}
        \begin{aligned}
                \mathbf{M}_{(\hat{\mathbf{i}},\hat{\mathbf{j}})} & = \int_{\left[0,1 \right]^3} \hat{\beta}_{\hat{\mathbf{i}}} \left(  \hat{x}  \right)  \, \hat{\beta}_{\hat{\mathbf{j}}} \left(  \hat{x}  \right)  \, \omega  \left(  \hat{x}  \right)  \, \mathrm{d} \hat{x}, \\
                \mathbf{K}_{(\hat{\mathbf{i}},\hat{\mathbf{j}})} &= \int_{\left[0,1 \right]^3}  \left(  Q  \left(  \hat{x}  \right)  \nabla \hat{\beta}_{\hat{\mathbf{i}}} \left(  \hat{x}  \right)  \right)  \cdot \nabla \hat{\beta}_{\hat{\mathbf{j}}} \left(  \hat{x}  \right)  \, \mathrm{d} \hat{x}\\
                & = \sum_{k,l=1}^{3} \int_{\left[0,1 \right]^3} q_{kl} \left( \hat{x} \right)  \frac{\partial}{\partial \hat{x}^{\left( l \right)}} \hat{\beta}_{\hat{\mathbf{i}}} \left( \hat{x} \right)  \frac{\partial}{\partial \hat{x}^{\left( k \right)}} \hat{\beta}_{\hat{\mathbf{j}}} \left( \hat{x} \right)  \, \mathrm{d} \hat{x},
        \end{aligned}
\end{equation}
for $\hat{\mathbf{i}}, \hat{\mathbf{j}} \in \hat{\mathbf{I}}$, where $\nabla$ denotes the gradient with respect to the parametric variables $\hat{x}$. In the following, we use the terms matrix and tensor equivalently, since they are related by reshaping.

The geometry-induced kernel functions, which we refer to as weight functions, are given by
{\small
\begin{equation}
        \label{eq:WeightFunctions}
        \omega \left( \hat{x} \right)  = \det \left( \nabla G \left( \hat{x} \right)  \right) , \quad Q \left( \hat{x} \right)  = \omega \left( \hat{x} \right)   \left( \nabla G \left( \hat{x} \right)  \right)^{-1}  \left( \nabla G \left( \hat{x} \right)  \right)^{-\top} =  \left(  q_{kl} \left( \hat{x} \right)   \right) _{k,l=1}^{3}.
\end{equation}}
The determinant is positive due to the orientation-preserving assumption \cref{eq:OrientationPreserving}. Although the basis functions $\hat{\beta}_{\hat{\mathbf{i}}}$ and their partial derivatives are separable, the geometry-induced weights are generally not. As shown in \cref{subsec:MassTensor}, for tensor-product B-spline geometries, $\omega$ is a spline function admitting a finite separated representation, whereas the entries $q_{kl}$ are generally rational and do not belong to a polynomial tensor-product spline space. Their nonseparability prevents a direct reduction of $\mathbf{K}_{(\hat{\mathbf{i}},\hat{\mathbf{j}})}$ to univariate integrals. Special geometries with affine parameterizations, such as the unit cube equipped with the identity map, constitute exceptions for which $\omega$ and the entries $q_{kl}$ of $Q$ are constant and therefore trivially separable. Consequently, when $\omega$ and $q_{kl}$ are treated as nonseparable functions, the integrands in \cref{eq:MassAndStiffness} are not separable, leading to computationally expensive multivariate quadrature in conventional assembly procedures.\\

To formulate the required separated approximations, we briefly introduce the low-rank tensor formats used throughout this work. A $D$-dimensional tensor $\mathbf{W}\in\mathbb{R}^{\left(n^{\left(1\right)},\ldots,n^{\left(D\right)}\right)}$ admits a canonical polyadic (CP) representation with $R$ terms if
\begin{equation}
        \label{eq:CPFormat}
        \mathbf{W} = \sum_{r=1}^{R} \bigotimes_{d=1}^{D} w_r^{\left(d\right)}, \qquad w_r^{\left(d\right)}\in\mathbb{R}^{n^{\left(d\right)}}.
\end{equation}
This notation directly expresses separability and the resulting Kronecker product structure. Each summand in \cref{eq:CPFormat} is a rank-one tensor and directly represents a product of univariate factors.

In this paper, however, the CP format is not used as the computational format, since computing CP approximations may be ill-posed and numerically unstable~\cite{tensorDecompositions}. Instead, all computations are performed in the TT format~\cite{Oseledets:2011}, defined by
\begin{equation}
        \label{eq:TT_format}
        \mathbf{W}_{\left(i^{\left(1\right)},\ldots,i^{\left(D\right)}\right)} = \mathbf{W}^{\left(1\right)}_{\left(:,i^{\left(1\right)},:\right)} \mathbf{W}^{\left(2\right)}_{\left(:,i^{\left(2\right)},:\right)} \cdots \mathbf{W}^{\left(D\right)}_{\left(:,i^{\left(D\right)},:\right)}.
\end{equation}
Here, $\mathbf{W}^{(d)}\in\mathbb{R}^{\left(R^{(d-1)},n^{(d)},R^{(d)} \right)}$ are the TT cores and the quantities\\ 
$R^{\left(1\right)},\ldots,R^{\left(D-1\right)}$ are the TT ranks with $R^{(0)}=R^{(D)}=1$. The storage requirement is $\sum_{d=1}^{D}R^{(d-1)}n^{(d)}R^{(d)}$, and TT rounding compresses a tensor in TT format by reducing its TT ranks while controlling the approximation error~\cite{Oseledets:2011}. Moreover, a TT tensor can be expanded into a CP representation with at most $\prod_{d=1}^{D-1} R^{\left(d\right)}$ rank-one terms. In this work we use the CP format \cref{eq:CPFormat} to display separability, while all computations are performed in the TT format. \\

The key idea, introduced in~\cite{angelos1}, to overcome the issue of the nonseparable integrands in \cref{eq:MassAndStiffness} is to approximate the non-separable weight functions by separable spline functions. For this purpose, let $\tilde{\mathbb{S}}$ be a three-variate tensor-product B-spline space with basis tensor
\begin{equation*}
        \tilde{\mathbf{B}} \left( \hat{x} \right)  = \bigotimes_{d=1}^{3} \tilde{b}^{\left( d \right)} \left( \hat{x}^{\left( d \right)} \right)  \in \mathbb{R}^{\left( \tilde{n}^{\left( 1 \right)},\tilde{n}^{\left( 2 \right)},\tilde{n}^{\left( 3 \right)} \right)},
\end{equation*}
where $\tilde{b}^{\left( d \right)} \left( \hat{x}^{\left( d \right)} \right)  = \left[ \tilde{\beta}^{\left( d \right)}_{1} \left(  \hat{x}^{\left( d \right)}  \right) , \ldots, \tilde{\beta}^{\left( d \right)}_{\tilde{n}^{\left( d \right)}} \left(  \hat{x}^{\left( d \right)}  \right)  \right]^\top$. The scalar weight $\omega$ is interpolated in this space by
\begin{equation*}
        \omega \left( \hat{x} \right)  \approx  \left\langle \mathbf{W}, \tilde{\mathbf{B}} \left( \hat{x} \right)   \right\rangle_F, \qquad \mathbf{W} \in \mathbb{R}^{\left( \tilde{n}^{\left( 1 \right)},\tilde{n}^{\left( 2 \right)},\tilde{n}^{\left( 3 \right)} \right)}.
\end{equation*}
The coefficient tensor $\mathbf{W}$ is determined by interpolation at the tensor-product grid of Greville points $\tilde{\mathbf{X}} = \tilde{X}^{\left( 1 \right)} \otimes \tilde{X}^{\left( 2 \right)} \otimes \tilde{X}^{\left( 3 \right)}$, cf.~\cite{greville}. Since both the basis $\tilde{\mathbf{B}}$ and the grid $\tilde{\mathbf{X}}$ are tensor products, interpolation leads to the following linear system 
\begin{equation}
        \label{eq:InterpolationSystem}
        \left(\tilde{b}^{\left( 1 \right)} \left( \tilde{X}^{\left( 1 \right)} \right) \otimes \tilde{b}^{\left( 2 \right)} \left( \tilde{X}^{\left( 2 \right)} \right) \otimes \tilde{b}^{\left( 3 \right)} \left( \tilde{X}^{\left( 3 \right)} \right) \right)^{\top} \ \operatorname{vec} \left( \mathbf{W} \right) = \operatorname{vec} \left( \omega \left( \mathbf{\tilde{X}} \right) \right).
\end{equation}
The full assembly of the linear system \cref{eq:InterpolationSystem} as well as the straightforward computation of $\mathbf{W}$ are in many situations too expensive. Following~\cite{BuengerDolgovStoll:2020}, we therefore never form the interpolation matrix explicitly and solve \cref{eq:InterpolationSystem} directly for a low-rank TT approximation $\tilde{\mathbf{W}}$ of $\mathbf{W}$ using \amen. We point out that assembling the right-hand side in \cref{eq:InterpolationSystem} requires evaluating the weight function $\omega$ at every point of the tensor grid $\tilde{\mathbf{X}}$, i.e., $\prod_{d=1}^{3} \tilde{n}^{\left( d \right)}$ evaluations, which becomes a computational bottleneck for fine discretizations.

The approximate solution of \amen~allows the separable representation
\begin{equation*}
        \mathbf{W} \approx \tilde{\mathbf{W}} = \sum_{r=1}^{R} \bigotimes_{d=1}^{3} w_r^{\left( d \right)}, \qquad w_r^{\left( d \right)}\in \mathbb{R}^{\tilde{n}^{\left( d \right)}},
\end{equation*}
so that the interpolant is separable
\begin{equation}
        \label{eq:SeparableWeight}
        \omega \left( \hat{x} \right)  \approx  \left\langle \tilde{\mathbf{W}}, \tilde{\mathbf{B}} \left( \hat{x} \right)   \right\rangle_F = \sum_{r=1}^{R} \prod_{d=1}^{3}  \left(  w_r^{\left( d \right)} \cdot \tilde{b}^{\left( d \right)} \left( \hat{x}^{\left( d \right)} \right)   \right) .
\end{equation}
Substituting \cref{eq:SeparableWeight} into the mass tensor in \cref{eq:MassAndStiffness} separates the three-dimensional integrals into products of univariate integrals. As a result, the mass tensor can be approximated by
\begin{equation}
        \label{eq:LowRankMassMatrix}
        \mathbf{M} \approx \sum_{r=1}^{R} \bigotimes_{d=1}^{3} M_r^{\left( d \right)},
\end{equation}
where
\begin{equation*}
        M_r^{\left( d \right)} = \int_{0}^{1} \hat{b}^{\left( d \right)} \left( \hat{x}^{\left( d \right)} \right)  \otimes \hat{b}^{\left( d \right)} \left( \hat{x}^{\left( d \right)} \right) \, \left(  w_r^{\left( d \right)} \cdot \tilde{b}^{\left( d \right)} \left( \hat{x}^{\left( d \right)} \right)   \right) \, \mathrm{d}\hat{x}^{\left( d \right)}.
\end{equation*}

The same procedure is applied to every entry of the matrix-valued weight $Q \left( \hat{x} \right)$. For $k,l=1,2,3$, we approximate
\begin{equation*}
        q_{kl} \left( \hat{x} \right)  \approx  \left\langle \tilde{\mathbf{V}}_{kl}, \tilde{\mathbf{B}} \left( \hat{x} \right)   \right\rangle_F = \sum_{r=1}^{R_{kl}} \prod_{d=1}^{3}  \left(  v_{kl,r}^{\left( d \right)} \cdot \tilde{b}^{\left( d \right)} \left( \hat{x}^{\left( d \right)} \right)   \right) ,
\end{equation*}
with $v_{kl,r}^{\left( d \right)}\in\mathbb{R}^{\tilde{n}^{\left( d \right)}}$. Introducing the univariate operator
\begin{equation}
        \label{eq:DifferentialStiffnessOperator}
        \delta \left( k,d \right)  f =
        \begin{cases}
                \dfrac{\partial f}{\partial \hat{x}^{\left( d \right)}}, & k=d,\\
                f, & k\neq d,
        \end{cases}
\end{equation}
the stiffness matrix can be approximated by
\begin{equation*}
        \mathbf{K} \approx \sum_{k,l=1}^{3} \sum_{r=1}^{R_{kl}} \bigotimes_{d=1}^{3} K_{kl,r}^{\left( d \right)},
\end{equation*}
where
{\small
\begin{equation*}
        K_{kl,r}^{\left( d \right)} = \int_{0}^{1} \left( \delta \left( l,d \right)  b^{\left( d \right)} \left( \hat{x}^{\left( d \right)} \right) \right)  \otimes \left( \delta \left( k,d \right)  b^{\left( d \right)} \left( \hat{x}^{\left( d \right)} \right) \right)  \, \left(  v_{kl,r}^{\left( d \right)} \cdot \tilde{b}^{\left( d \right)} \left( \hat{x}^{\left( d \right)} \right)   \right) \, \mathrm{d}\hat{x}^{\left( d \right)}.
\end{equation*}}
This means that after the interpolation of the non-separable weights, the assembly of the multivariate mass and stiffness matrices in \cref{eq:MassAndStiffness} is reduced to the assembly of small univariate matrices.

For both approximations the interpolation space $\tilde{\mathbb{S}}$ must be rich enough to represent the weight functions with the desired accuracy. For the mass matrix, an interpolation space that contains $\omega$ can be derived from the geometry space itself. As shown in \cref{subsec:MassTensor}, each summand of $\omega$ is, in every parametric direction $d$, the product of two degree-$p^{(d)}$ splines and one degree-$(p^{(d)}-1)$ derivative, i.e.\ a piecewise polynomial of degree $3p^{(d)}-1$. By the same regularity argument used in \cref{subsec:MassTensor}, all such products lie in a spline space $\mathbb{S}_{3p^{(d)}-1}\left(\tilde{\Xi}^{(d)}\right)$, whose auxiliary knot vector $\tilde{\Xi}^{(d)}$ has interior multiplicities $2p^{(d)}+\mu^{(d)}_\ell$ and boundary multiplicities $3p^{(d)}$. As a consequence, $\omega$ belongs to the tensor-product spline space
\begin{equation}
        \label{eq:WeightFunctionSplineSpace}
        \mathbb{S}_{3p^{(1)}-1} \left(\tilde{\Xi}^{(1)}\right) \otimes \mathbb{S}_{3p^{(2)}-1} \left(\tilde{\Xi}^{(2)}\right) \otimes \mathbb{S}_{3p^{(3)}-1} \left(\tilde{\Xi}^{(3)}\right),
\end{equation}
and can be recovered exactly by interpolation in this space, provided the interpolation system is nonsingular. Any remaining error is then due only to the low-rank truncation, numerical errors in the solution process of the interpolation system or the use of a smaller auxiliary space. 

In contrast, the rational functions $q_{kl}$ generally do not belong to any polynomial spline space. Following~\cite{angelos1}, one may start from the knot vectors used for $\omega$ and increase the spline degrees until the desired accuracy is reached. The required degree depends on the geometry and the prescribed accuracy, and increasing it must be balanced against the resulting growth of the interpolation system via the Greville points and the coefficient tensor, which constitutes the main computational limitation of the interpolation approach.

\section{Projection-based low-rank method}
\label{sec:Projection}

In this section, we present the proposed low-rank assembly methods for the mass tensor $\mathbf{M}$ and the stiffness tensor $\mathbf{K}$. Both exploit the tensor-product structure of the B-spline geometry map to rewrite $\omega$ and $q_{kl}$, $k,l=1,2,3$, of the matrix-valued weight function $Q$ in separated form, reducing the three-variate Galerkin integrals \cref{eq:MassAndStiffness} to products of univariate integrals compatible with the TT format. The common mechanism is that products of univariate B-splines can again be represented exactly in univariate spline spaces of higher degree, where in both cases the corresponding coefficient transfer matrices are obtained by univariate $L^2$-projections. For the mass tensor in \cref{subsec:MassTensor}, the weight $\omega$ is piecewise polynomial, so this yields an exact assembly, up to TT truncation and quadrature. For the stiffness tensor in \cref{subsec:StiffTensor}, each entry $q_{kl}$ splits into a polynomial numerator and the common reciprocal determinant, which is rational and has to be approximated by a projection.

\subsection{Low-rank mass tensor}
\label{subsec:MassTensor}
In the following, let $G \colon [0,1]^3 \to \Omega$ be the geometry map representing the computational domain $\Omega \subset \mathbb{R}^3$, constructed from the tensor-product B-spline space $\mathbb{S}_{p} \left( \Xi^{\left( 1 \right)},\Xi^{\left( 2 \right)},\Xi^{\left( 3 \right)} \right)$. The aim is to bring the mass tensor \\
$\mathbf{M} \in \mathbb{R}^{\left( \hat{n}^{(1)}, \hat{n}^{(2)}, \hat{n}^{(3)} \right) \times \left( \hat{n}^{(1)}, \hat{n}^{(2)}, \hat{n}^{(3)} \right)}$ defined by \cref{eq:MassAndStiffness} into a separable form. Using the Leibniz formula, the determinant can be written as
\begin{equation}
        \label{eq:WeightFunctionDeterminant}
        \begin{aligned}
                \omega \left( \hat{x} \right) &= \det \left( \nabla G \left( \hat{x} \right)  \right)  \\
                & = \frac{\partial G_1}{\partial \hat{x}^{\left( 1 \right)}} \left( \hat{x} \right) \frac{\partial G_2}{\partial \hat{x}^{\left( 2 \right)}} \left( \hat{x} \right) \frac{\partial G_3}{\partial \hat{x}^{\left( 3 \right)}} \left( \hat{x} \right) + \frac{\partial G_1}{\partial \hat{x}^{\left( 2 \right)}} \left( \hat{x} \right) \frac{\partial G_2}{\partial \hat{x}^{\left( 3 \right)}} \left( \hat{x} \right) \frac{\partial G_3}{\partial \hat{x}^{\left( 1 \right)}} \left( \hat{x} \right)  \\
                & \quad + \frac{\partial G_1}{\partial \hat{x}^{\left( 3 \right)}} \left( \hat{x} \right) \frac{\partial G_2}{\partial \hat{x}^{\left( 1 \right)}} \left( \hat{x} \right)  \frac{\partial G_3}{\partial \hat{x}^{\left( 2 \right)}} \left( \hat{x} \right) - \frac{\partial G_1}{\partial \hat{x}^{\left( 3 \right)}} \left( \hat{x} \right)  \frac{\partial G_2}{\partial \hat{x}^{\left( 2 \right)}} \left( \hat{x} \right)  \frac{\partial G_3}{\partial \hat{x}^{\left( 1 \right)}} \left( \hat{x} \right)  \\
                & \quad - \frac{\partial G_1}{\partial \hat{x}^{\left( 1 \right)}} \left( \hat{x} \right) \frac{\partial G_2}{\partial \hat{x}^{\left( 3 \right)}} \left( \hat{x} \right) \frac{\partial G_3}{\partial \hat{x}^{\left( 2 \right)}} \left( \hat{x} \right)  - \frac{\partial G_1}{\partial \hat{x}^{\left( 2 \right)}} \left( \hat{x} \right) \frac{\partial G_2}{\partial \hat{x}^{\left( 1 \right)}} \left( \hat{x} \right) \frac{\partial G_3}{\partial \hat{x}^{\left( 3 \right)}} \left( \hat{x} \right) .
        \end{aligned}
\end{equation}
Since the components of the geometry map are Frobenius products of the B-spline basis tensor and the coordinate control-point tensors, see \cref{eq:GeometryMapTensor}, we obtain, for example,
\begin{equation}
        \label{eq:GeometryMapPartialDerivative}
        \frac{\partial G_1}{\partial \hat{x}^{\left( 2 \right)}} \left( \hat{x} \right) = \frac{\partial}{\partial \hat{x}^{\left( 2 \right)}} \left\langle \mathbf{C}^{(1)}, \mathbf{B} \left( \hat{x} \right) \right\rangle_F = \left\langle \mathbf{C}^{(1)}, \partial^{\left( 2 \right)} \mathbf{B} \left( \hat{x} \right) \right\rangle_F .
\end{equation}
Here, $\mathbf{C}^{(a)} = \mathbf{C}_{\left( a,:,:,: \right)} \in \mathbb{R}^{\left( n^{\left( 1 \right)}, n^{\left( 2 \right)}, n^{\left( 3 \right)} \right)}$, $a=1,2,3$, and we introduce the notation
\begin{equation*}
        \partial^{\left( 2 \right)} \mathbf{B} \left( \hat{x} \right) = b^{\left( 1 \right)} \left( \hat{x}^{\left( 1 \right)} \right) \otimes \frac{\partial}{\partial \hat{x}^{\left( 2 \right)}} b^{\left( 2 \right)} \left( \hat{x}^{\left( 2 \right)} \right) \otimes b^{\left( 3 \right)} \left( \hat{x}^{\left( 3 \right)} \right) .
\end{equation*}
The partial derivatives $\partial^{\left( 1 \right)} \mathbf{B} \left( \hat{x} \right)$ and $\partial^{\left( 3 \right)} \mathbf{B} \left( \hat{x} \right)$ are defined analogously. With this notation, the first summand in \cref{eq:WeightFunctionDeterminant}, for example, can be expressed as
{\small
\begin{align*}
        \frac{\partial G_1}{\partial \hat{x}^{\left( 1 \right)}} \left( \hat{x} \right) \frac{\partial G_2}{\partial \hat{x}^{\left( 2 \right)}} \left( \hat{x} \right) \frac{\partial G_3}{\partial \hat{x}^{\left( 3 \right)}} \left( \hat{x} \right)  & = \left\langle \mathbf{C}^{(1)}, \partial^{\left( 1 \right)} \mathbf{B} \left( \hat{x} \right)  \right\rangle_F \left\langle \mathbf{C}^{(2)}, \partial^{\left( 2 \right)} \mathbf{B} \left( \hat{x} \right) \right\rangle_F \left\langle \mathbf{C}^{(3)}, \partial^{\left( 3 \right)} \mathbf{B} \left( \hat{x} \right) \right\rangle_F \\
        & = \left\langle \mathbf{C}^{(1)} \otimes \mathbf{C}^{(2)} \otimes \mathbf{C}^{(3)}, \partial^{\left( 1 \right)} \mathbf{B} \left( \hat{x} \right) \otimes \partial^{\left( 2 \right)} \mathbf{B} \left( \hat{x} \right) \otimes \partial^{\left( 3 \right)} \mathbf{B} \left( \hat{x} \right) \right\rangle_F 
\end{align*}}
where we use the fact that a product of Frobenius products can be written as one Frobenius product involving the tensor products of the corresponding arguments.

By using the bilinearity of the Frobenius product, we can rewrite \cref{eq:WeightFunctionDeterminant} as
\begin{equation}
        \label{eq:WeightFunctionFrobeniusProduct}
        \omega \left( \hat{x} \right) = \left\langle \mathbf{C}_{\Sigma}, \partial^{\left( 1 \right)} \mathbf{B} \left( \hat{x} \right) \otimes \partial^{\left( 2 \right)} \mathbf{B} \left( \hat{x} \right) \otimes \partial^{\left( 3 \right)} \mathbf{B} \left( \hat{x} \right) \right\rangle_F ,
\end{equation}
where
\begin{equation}
        \label{eq:WeightFunctionFrobeniusProductCoefficients}
        \begin{aligned}
                \mathbf{C}_{\Sigma} &= \mathbf{C}^{(1)} \otimes \mathbf{C}^{(2)} \otimes \mathbf{C}^{(3)} + \mathbf{C}^{(3)} \otimes \mathbf{C}^{(1)} \otimes \mathbf{C}^{(2)} + \mathbf{C}^{(2)} \otimes \mathbf{C}^{(3)} \otimes \mathbf{C}^{(1)} \\
                & \quad - \mathbf{C}^{(3)} \otimes \mathbf{C}^{(2)} \otimes \mathbf{C}^{(1)} - \mathbf{C}^{(1)} \otimes \mathbf{C}^{(3)} \otimes \mathbf{C}^{(2)} - \mathbf{C}^{(2)} \otimes \mathbf{C}^{(1)} \otimes \mathbf{C}^{(3)} .
        \end{aligned}
\end{equation}
As a result,
\begin{equation*}
        \mathbf{C}_{\Sigma} \in \mathbb{R}^{\left( n^{\left( 1 \right)}, n^{\left( 2 \right)}, n^{\left( 3 \right)}, n^{\left( 1 \right)}, n^{\left( 2 \right)}, n^{\left( 3 \right)}, n^{\left( 1 \right)}, n^{\left( 2 \right)}, n^{\left( 3 \right)} \right)}.
\end{equation*}
Inserting this representation into the definition of the mass tensor in \cref{eq:MassAndStiffness} yields
{\small
\begin{equation}
        \label{eq:MassMatrixLowrankPermutation}
        \begin{aligned}
                \mathbf{M} &= \int_{[0,1]^3} \hat{\mathbf{B}} \left( \hat{x} \right) \otimes \hat{\mathbf{B}} \left( \hat{x} \right) \, \left\langle \mathbf{C}_{\Sigma}, \partial^{\left( 1 \right)} \mathbf{B} \left( \hat{x} \right) \otimes \partial^{\left( 2 \right)} \mathbf{B} \left( \hat{x} \right) \otimes \partial^{\left( 3 \right)} \mathbf{B} \left( \hat{x} \right) \right\rangle_F \,\mathrm{d}\hat{x} \\
                & = \int_{[0,1]^3} \left( \hat{\mathbf{B}} \left( \hat{x} \right) \otimes \hat{\mathbf{B}} \left( \hat{x} \right) \otimes \partial^{\left( 1 \right)} \mathbf{B} \left( \hat{x} \right) \otimes \partial^{\left( 2 \right)} \mathbf{B} \left( \hat{x} \right)  \otimes \partial^{\left( 3 \right)} \mathbf{B} \left( \hat{x} \right) \right) \cdot \mathbf{C}_{\Sigma} \,\mathrm{d}\hat{x} \\
                & = \int_{[0,1]^3} \left( \mathbf{B}^{\left( 1 \right)} \left( \hat{x}^{\left( 1 \right)} \right) \otimes \mathbf{B}^{\left( 2 \right)} \left( \hat{x}^{\left( 2 \right)} \right) \otimes \mathbf{B}^{\left( 3 \right)} \left( \hat{x}^{\left( 3 \right)} \right) \right) \cdot \mathbf{C}_{\Sigma} \,\mathrm{d}\hat{x}.
        \end{aligned}
\end{equation}}
Here, $\hat{\mathbf{B}} \left( \hat{x} \right)$ denotes the tensor containing all basis functions of the solution space $V_h$, which has a tensor-product structure. In the second line of \cref{eq:MassMatrixLowrankPermutation}, we use that the product of $\hat{\mathbf{B}} \left( \hat{x} \right)  \otimes \hat{\mathbf{B}} \left( \hat{x} \right)$ with the Frobenius product can be written as a contracted product, denoted by $\cdot$. In the last line, the tensor factors are permuted and grouped according to the dimensions and $\mathbf{C}_{\Sigma}$ is permuted accordingly. We define
{\small
\begin{equation}
        \label{eq:BsplineTensorBeforeProjection}
        \mathbf{B}^{\left( d \right)} \left( \hat{x}^{\left( d \right)} \right) = \hat{b}^{\left( d \right)} \left( \hat{x}^{\left( d \right)} \right) \otimes \hat{b}^{\left( d \right)} \left( \hat{x}^{\left( d \right)} \right) \otimes b^{\left( d \right)} \left( \hat{x}^{\left( d \right)} \right) \otimes b^{\left( d \right)} \left( \hat{x}^{\left( d \right)} \right) \otimes \frac{\partial}{\partial \hat{x}^{\left( d \right)}} b^{\left( d \right)} \left( \hat{x}^{\left( d \right)} \right) ,
\end{equation}}
for $d=1,2,3$. We note that, $\mathbf{B}^{\left( d \right)} \left( \hat{x}^{\left( d \right)} \right)  \in \mathbb{R}^{\left( \hat{n}^{\left( d \right)},\hat{n}^{\left( d \right)},n^{\left( d \right)},n^{\left( d \right)},n^{\left( d \right)} \right)}$. Here, $\hat{n}^{\left( d \right)}$ denotes the number of univariate basis functions in dimension $d$ of the solution space $V_h$, while $n^{\left( d \right)}$ denotes the corresponding number of univariate basis functions used for the geometry map.

Using a low-rank approximation of the coefficient tensor
\begin{equation}
        \label{eq:WeightFunctionFrobeniusProductCoefficientsLowRank}
        \mathbf{C}_{\Sigma} \approx \tilde{\mathbf{C}}_{\Sigma} = \sum_{r=1}^{R_{\Sigma}} \bigotimes_{d=1}^{3} \tilde{\mathbf{C}}^{(d)}_{\Sigma,r}, \qquad \tilde{\mathbf{C}}^{(d)}_{\Sigma,r} \in \mathbb{R}^{\left( n^{\left( d \right)},n^{\left( d \right)},n^{\left( d \right)} \right)},
\end{equation}
the mass tensor can be represented approximately as a sum of tensor products,
{\small
\begin{equation}
        \label{eq:MassMatrixLowrank_1}
        \mathbf{M} \approx \sum_{r=1}^{R_{\Sigma}} \bigotimes_{d=1}^{3} \left( \tilde{\mathbf{C}}^{(d)}_{\Sigma,r} \cdot  \int_0^1 \mathbf{B}^{\left( d \right)} \left( \hat{x}^{\left( d \right)} \right) \,\mathrm{d}\hat{x}^{\left( d \right)} \right) \in \mathbb{R}^{\left( \hat{n}^{(1)}, \hat{n}^{(2)}, \hat{n}^{(3)} \right) \times \left( \hat{n}^{(1)}, \hat{n}^{(2)}, \hat{n}^{(3)} \right)}
\end{equation}}
where the contracted product $\cdot$ acts on the last three modes of $\mathbf{B}^{\left( d \right)} \left( \hat{x}^{\left( d \right)} \right)$ and on all three modes of $\tilde{\mathbf{C}}^{(d)}_{\Sigma,r}$. Therefore, each factor in \cref{eq:MassMatrixLowrank_1} is a matrix in $\mathbb{R}^{\hat{n}^{\left( d \right)} \times \hat{n}^{\left( d \right)}}$. \\

Although \cref{eq:MassMatrixLowrank_1} has a separated structure, this formulation leads to two challenges. First, looking at \cref{eq:WeightFunctionFrobeniusProductCoefficients}, one would have to compute a low-rank representation of a sum of tensors of order $9$, which can be computationally challenging. Second, the integration of $\mathbf{B}^{\left( d \right)} \left( \hat{x}^{\left( d \right)} \right)$, which is a tensor of order five, remains a computational bottleneck. Both difficulties are caused by the modes associated with the basis of the geometry map $G$. More precisely, the first two modes in \cref{eq:BsplineTensorBeforeProjection} are given by $\hat{b}^{\left( d \right)} \left( \hat{x}^{\left( d \right)} \right)$ and correspond to the solution space $V_h$. Since these modes determine the richness of the discrete solution space, they should be preserved. The remaining three modes, however, are induced by the geometry representation and have size $\left( n^{\left( d \right)} \right)^3$.

As a consequence, for each pair of univariate basis functions from the solution space, the formulation still involves $\left( n^{\left( d \right)} \right)^3$ products of geometry basis functions in dimension $d$. In addition, for a rank-$R$ separated representation, one has to work with $R \left( n^{\left( d \right)} \right)^3$ many corresponding coefficients. Therefore, the aim is to compress the last three modes of \cref{eq:BsplineTensorBeforeProjection}, while leaving the two modes associated with the solution space $V_h$ unchanged.

The last three modes of \cref{eq:BsplineTensorBeforeProjection} are products of univariate B-splines from $\mathbb{S}_{p^{\left( d \right)}} \left( \Xi^{\left( d \right)} \right)$ and contain $\left( n^{\left( d \right)} \right)^3$ entries. A single entry is given by
{\small
\begin{equation}
        \label{eq:SplineProduct3p-1}
        \begin{gathered}
                \left(  b^{\left( d \right)} \left( \hat{x}^{\left( d \right)} \right)  \otimes b^{\left( d \right)} \left( \hat{x}^{\left( d \right)} \right) \otimes \frac{\partial}{\partial \hat{x}^{\left( d \right)}} b^{\left( d \right)} \left( \hat{x}^{\left( d \right)} \right) \right) _{\left( i,j,k \right)}  \\
                = \beta_i^{\left( d \right)} \left( \hat{x}^{\left( d \right)} \right) \beta_j^{\left( d \right)} \left( \hat{x}^{\left( d \right)} \right) \left( \beta_k^{\left( d \right)} \right)' \left( \hat{x}^{\left( d \right)} \right),
        \end{gathered}
\end{equation}}
which is the product of two splines of degree $p^{\left( d \right)}$ and one derivative of a spline of degree $p$. Since the derivative can be represented in a spline space of degree $p^{\left( d \right)}-1$, see \cref{eq:BSplineDerivative}, the product in \cref{eq:SplineProduct3p-1} can be represented in a B-spline basis of degree $3p^{\left( d \right)}-1$ and reduced regularity, as already known from~\cite{angelos1}. \\

Let
\begin{equation*}
        \Xi^{(d)} = \left\{ \xi^{(d)}_1,\ldots, \xi^{(d)}_{n^{(d)}+p^{(d)}+1} \right\}
\end{equation*}
be the open knot vector of the univariate tensor-product geometry space $\mathbb{S}_{p^{(d)}}\left(\Xi^{(d)}\right)$. We denote the distinct knot values, often referred to as breaking points, by $\zeta^{(d)}_0 < \zeta^{(d)}_1 < \cdots < \zeta^{(d)}_{m^{(d)}}$, and their multiplicities in $\Xi^{(d)}$ by $\mu^{(d)}_\ell$, $\ell=0,\ldots,m^{(d)}$. The boundary knots are $\zeta^{(d)}_0 = 0$ and $\zeta^{(d)}_{m^{(d)}} = 1$ with multiplicity $\mu^{(d)}_0 = \mu^{(d)}_{m^{(d)}} = p^{(d)}+1$, while the distinct interior knots are $\zeta^{(d)}_1,\ldots,\zeta^{(d)}_{m^{(d)}-1}$. The derivative of a spline in $\mathbb{S}_{p^{(d)}}\left(\Xi^{(d)}\right)$ belongs to the spline space $\mathbb{S}_{p^{(d)}-1}\left(\Xi_{\partial}^{(d)}\right)$ with knot vector
\begin{equation}
        \label{eq:DerivativeKnotvector}
        \Xi_{\partial}^{(d)} = \left\{ \xi^{(d)}_2,\ldots, \xi^{(d)}_{n^{(d)}+p^{(d)}} \right\}.
\end{equation}
As a result, one knot is removed at each endpoint, while the interior knot multiplicities remain unchanged. This is consistent with the regularity of the derivative. Indeed, at an interior knot $\zeta^{(d)}_\ell$ with multiplicity $\mu^{(d)}_\ell$, a spline of degree $p^{(d)}$ has regularity $C^{p^{(d)}-\mu^{(d)}_\ell}$. After differentiation, the regularity is reduced by one and becomes $C^{p^{(d)}-\mu^{(d)}_\ell-1}$. 

The product of two splines from $\mathbb{S}_{p^{(d)}}(\Xi^{(d)})$ is a spline of degree $2p^{(d)}$. At an interior knot $\zeta^{(d)}_\ell$ with multiplicity $\mu^{(d)}_\ell$, a spline of degree $p^{(d)}$ has continuity $C^{p^{(d)}-\mu^{(d)}_\ell}$ and therefore the product of two such splines has the same regularity. Thus, it can be represented in a degree-$2p^{(d)}$ spline space whose interior knot multiplicity $\mu^{(d)}_{\ell,2}$ satisfies
\begin{equation*}
        \mu^{(d)}_{\ell,2} = \mu^{(d)}_\ell + p^{(d)}, \qquad \ell=1,\ldots,m^{(d)}-1.
\end{equation*}
At the boundary knots, we use the open endpoint multiplicity $2p^{(d)}+1$ and the knot vector for the product space $\mathbb{S}_{2 p^{\left( d \right)}} \left( \Xi^{(d)}_2 \right)$ is given by
{\small
\begin{equation}
        \label{eq:ProductKnotvector}
        \Xi^{(d)}_2 = \big\{ \underbrace{\zeta^{(d)}_0,\ldots,\zeta^{(d)}_0}_{2p^{(d)}+1}, \ldots, \underbrace{\zeta^{(d)}_\ell,\ldots,\zeta^{(d)}_\ell}_{p^{(d)}+\mu^{(d)}_\ell}, \ldots, \underbrace{\zeta^{(d)}_{m^{(d)}},\ldots,\zeta^{(d)}_{m^{(d)}}}_{2p^{(d)}+1}\big\}.
\end{equation}}

We now multiply a function from this product space $\mathbb{S}_{2 p^{\left( d \right)}} \left( \Xi^{(d)}_2 \right)$ with a derivative basis function from $\mathbb{S}_{p^{\left( d \right)} - 1} \left( \Xi_{\partial}^{(d)} \right)$. The first factor has degree $2p^{(d)}$, while the derivative factor has degree $p^{(d)}-1$. Therefore, the resulting product has degree $3p^{(d)}-1$. The regularity of the product is determined by the less regular factor. The product spline has continuity $C^{p^{(d)}-\mu^{(d)}_\ell}$ at an interior knot, while the derivative spline has continuity $C^{p^{(d)}-\mu^{(d)}_\ell-1}$. This means that the triple product has also continuity $C^{p^{(d)}-\mu^{(d)}_\ell-1}$. For a spline of degree $3p^{(d)}-1$, this continuity corresponds to the interior knot multiplicity $\tilde{\mu}^{(d)}_\ell$ satisfying
\begin{equation*}
        \tilde{\mu}^{(d)}_\ell = 2p^{(d)} + \mu^{(d)}_\ell, \qquad \ell=1,\ldots,m^{(d)}-1.
\end{equation*}
At the endpoints, the open knot multiplicity for degree $3p^{(d)}-1$ is $3p^{(d)}$. Therefore, the reduced knot vector is
{\small
\begin{equation}
        \label{eq:ReducedTripleProductKnotVector}
        \tilde{\Xi}^{(d)} = \big\{ \underbrace{\zeta^{(d)}_0,\ldots,\zeta^{(d)}_0}_{3p^{(d)}}, \ldots, \underbrace{\zeta^{(d)}_\ell, \ldots, \zeta^{(d)}_\ell}_{2p^{(d)}+\mu^{(d)}_\ell}, \ldots, \underbrace{\zeta^{(d)}_{m^{(d)}},\ldots,\zeta^{(d)}_{m^{(d)}}}_{3p^{(d)}} \big\}.
\end{equation}}
By construction, the space $\mathbb{S}_{3p^{(d)}-1} \left(\tilde{\Xi}^{(d)}\right)$ contains all products \cref{eq:SplineProduct3p-1}. Therefore the weight function $\omega$ lies in the tensor-product spline space already introduced in \cref{eq:WeightFunctionSplineSpace}. For fixed degree $p^{(d)}$, the number of basis functions in this reduced space grows only linearly in $n^{(d)}$, in contrast to the cubic growth of the triple product in \cref{eq:SplineProduct3p-1}.\\

The next step is to bring the coefficient tensor $\tilde{\mathbf{C}}_{\Sigma}$, which is of order $9$, see \cref{eq:WeightFunctionFrobeniusProductCoefficientsLowRank}, into a grouped tensor of order $3$, whose modes correspond to the three parametric directions and later serve as input for the transfer into the reduced spline space \cref{eq:WeightFunctionSplineSpace}. The first challenge here is to construct $\tilde{\mathbf{C}}_{\Sigma}$ itself. In our method, this tensor is constructed in tensor train format, which is advantageous because the TT format behaves differently from the CP format when tensor products are formed.

We decompose the coordinate control-point tensors $\mathbf{C}^{(a)}$, $a = 1, 2, 3$, into TT format, form the tensor products and signed sum in \cref{eq:WeightFunctionFrobeniusProductCoefficients} directly in TT format, and finally permute and reshape the result so that the three modes belonging to the same parametric direction are grouped together. Unlike the CP format, where the tensor product of two rank-$R_1$ and rank-$R_2$ tensors has CP rank $R_1R_2$, tensor products in TT format are obtained by simply concatenating the TT cores of the two factors, so that ranks add rather than multiply. More precisely, if two tensors in TT format have storage complexities $\mathcal{O} \left( D R_1^2 n_{\max} \right)$ and $\mathcal{O} \left( D R_2^2 n_{\max} \right)$, then their tensor product can be represented by simply appending the TT cores of the second tensor to those of the first, and, before rounding, the storage complexity of the result is the sum of the two individual complexities, $\mathcal{O} \left( D R_1^2 n_{\max} + D R_2^2 n_{\max} \right)$, where $n_{\max} = \max_{d = 1, \ldots, D} n^{(d)}$ and $R_1, R_2$ are the maximal TT rank of the corresponding tensor.

In contrast to the CP format, however, a permutation of the modes can increase the TT ranks of a tensor. It is therefore advisable to round the tensor again after the permutation, as we do in the implementation. In the final step, we reshape the tensor so that all quantities depending on the same dimension $d$ are collected in one mode. The grouped coefficient tensor is thus in TT format with cores of size $\left( R^{\left( d-1 \right)}_{\Sigma}, \left( n^{\left( d \right)} \right)^3, R^{\left( d \right)}_{\Sigma} \right)$, $d=1,2,3$. Since a TT core has a single physical index, the mode of size $\left(n^{(d)}\right)^3$ is the multi-index $(i,j,k)$ of the triple product \cref{eq:SplineProduct3p-1}, and we address the entries accordingly. The resulting TT ranks originate from the TT ranks of the coordinate control-point tensors $\mathbf{C}^{(a)}$ and are mainly increased by the additions in \cref{eq:WeightFunctionFrobeniusProductCoefficients} and by the subsequent permutation. We point out that TT rounding only acts on the ranks between the three grouped modes. Since the indices belonging to one parametric direction are collected into a single physical mode of size $\left(n^{(d)}\right)^3$, rounding cannot reduce this mode size itself, only the TT ranks $R_\Sigma^{(d-1)}, R_\Sigma^{(d)}$ coupling it to the neighboring modes. To keep the notation simple, we write the resulting grouped tensor in CP-like form,
\begin{equation}
        \label{eq:WeightFunctionFrobeniusProductCoefficientsTT}
        \tilde{\mathbf{C}}_{\Sigma} = \sum_{r = 1}^{R_{\Sigma}} \bigotimes^{3}_{d = 1} \tilde{C}^{(d)}_{\Sigma, r}, \qquad \tilde{C}^{(d)}_{\Sigma, r} \in \mathbb{R}^{\left( n^{\left( d \right)} \right)^3},
\end{equation}
while all computations are carried out in TT format. \\

Next, we describe how the coefficients $\tilde{\mathbf{C}}_{\Sigma}$ are transferred to the spline space $\mathbb{S}_{3p^{(d)}-1}\left(\tilde{\Xi}^{(d)}\right)$. For this purpose, we construct, for each $d=1,2,3$, univariate coefficient transfer operators which map coefficients with respect to the products \cref{eq:SplineProduct3p-1} to coefficients with respect to the basis of $\mathbb{S}_{3p^{(d)}-1}\left(\tilde{\Xi}^{(d)}\right)$. 

Using the derivative, product, and the reduced knot vectors \cref{eq:DerivativeKnotvector,eq:ProductKnotvector,eq:ReducedTripleProductKnotVector} introduced above for the spaces $\mathbb{S}_{p^{(d)}-1}\left(\Xi_{\partial}^{(d)}\right)$, $\mathbb{S}_{2 p^{\left( d \right)}} \left( \Xi^{(d)}_2 \right)$ and $\mathbb{S}_{3p^{(d)}-1} \left(\tilde{\Xi}^{(d)}\right)$, and denoting the corresponding bases by $\left\{\beta_{\partial,b}^{(d)} \right\}_{b=1}^{n_\partial^{(d)}}$, $\left\{\beta_{2,a}^{(d)} \right\}_{a=1}^{n_2^{(d)}}$, and $\left\{\tilde{\beta}_\alpha^{(d)} \right\}_{\alpha=1}^{\tilde{n}^{(d)}}$, then one has the following representations
\begin{equation}
        \label{eq:CoefficientTransferMatrices}
        \begin{gathered}
                \left(\beta_k^{(d)}\right)' = \sum_{b=1}^{n_\partial^{(d)}} \Delta_{b,k}^{(d)}\, \beta_{\partial,b}^{(d)} \in \mathbb{S}_{p^{(d)}-1}\left(\Xi_{\partial}^{(d)}\right), \\
                \beta_i^{(d)} \beta_j^{(d)} = \sum_{a=1}^{n_2^{(d)}} P_{a,(i,j)}^{(d)}\, \beta_{2,a}^{(d)} \in \mathbb{S}_{2 p^{\left( d \right)}} \left( \Xi^{(d)}_2 \right), \\
                \beta_{2,a}^{(d)} \beta_{\partial,b}^{(d)} = \sum_{\alpha=1}^{\tilde{n}^{(d)}} \left(T_{\mathrm{mix}}^{(d)}\right)_{\alpha,(a,b)} \tilde{\beta}_\alpha^{(d)} \in \mathbb{S}_{3p^{(d)}-1} \left(\tilde{\Xi}^{(d)}\right).
        \end{gathered}
\end{equation}
Here, $\Delta^{(d)}\in\mathbb{R}^{n_\partial^{(d)}\times n^{(d)}}$ follows directly from \cref{eq:BSplineDerivative} and is exact by construction, while $P^{(d)}\in\mathbb{R}^{n_2^{(d)}\times n_{\mathrm{pair}}^{(d)}}$ and $T_{\mathrm{mix}}^{(d)}\in\mathbb{R}^{\tilde{n}^{(d)}\times n_{\mathrm{mix}}^{(d)}}$ are obtained by $L^2$-projection onto $\mathbb{S}_{2p^{(d)}}\left(\Xi_2^{(d)}\right)$ and $\mathbb{S}_{3p^{(d)}-1}\left(\tilde{\Xi}^{(d)}\right)$, respectively, and are exact up to quadrature since the target spaces contain the corresponding products by construction. In both cases, only overlapping index pairs $(i,j)$ and $(a,b)$ are considered, since the products vanish otherwise, and $n_{\mathrm{pair}}^{(d)}$ and $n_{\mathrm{mix}}^{(d)}$ denote the resulting numbers of overlapping pairs. This restriction to overlapping pairs is natural and avoids coefficient computations for products which vanish, since $\beta_i^{(d)}\beta_j^{(d)}=0$ whenever the supports of $\beta_i^{(d)}$ and $\beta_j^{(d)}$ do not overlap.

Composing the three maps in \cref{eq:CoefficientTransferMatrices} gives the coefficient transfer from the triple product \cref{eq:SplineProduct3p-1} to the basis of $\mathbb{S}_{3p^{(d)}-1}\left(\tilde{\Xi}^{(d)}\right)$. Rather than assembling this composition as one dense matrix, we apply $P^{(d)}$, $\Delta^{(d)}$, and $T_{\mathrm{mix}}^{(d)}$ successively and directly to each $\tilde{C}^{(d)}_{\Sigma,r} \in \mathbb{R}^{\left( n^{\left( d \right)} \right)^3}$ in \cref{eq:WeightFunctionFrobeniusProductCoefficientsTT}. Rather than assembling this composition as one dense matrix, we apply $P^{(d)}$, $\Delta^{(d)}$, and $T_{\mathrm{mix}}^{(d)}$ successively and directly to each $\tilde{\mathbf{C}}^{(d)}_{\Sigma,r}$ in \cref{eq:WeightFunctionFrobeniusProductCoefficientsTT}. Since $P^{(d)}$ and $\Delta^{(d)}$ act separately on the pair $(i,j)$ and on $k$, we write $X = \tilde{\mathbf{C}}^{(d)}_{\Sigma,r}$ with entries $X_{(i,j,k)}$ and compute
\begin{equation}
        \label{eq:SequentialPushforward}
        \begin{gathered}
                Y_{(a,k)} = \sum_{i,j} P^{(d)}_{a,(i,j)}\, X_{(i,j,k)}, \\
                Z_{(a,\mu)} = \sum_{k} \Delta^{(d)}_{\mu,k}\, Y_{(a,k)}, \\
                g_{r,\alpha}^{(d)} = \sum_{(a,\mu)} \left(T_{\mathrm{mix}}^{(d)}\right)_{\alpha,(a,\mu)} Z_{(a,\mu)}.
        \end{gathered}
\end{equation}
The computation of $Y$ sums over overlapping pairs $(i,j)$ and can be understood as a contracted product acting on the first two dimensions of $X$, with the contraction coefficients stored in the columns of $P^{(d)}$, restricted to overlapping pairs. Since $Y$ still depends on both $a$ and $k$, it is a matrix rather than a vector, and the computation of $Z$ is accordingly a genuine matrix-matrix product, $Z^{(d)} = Y^{(d)} \left(\Delta^{(d)}\right)^\top$. The computation of $g_r^{(d)}$ is a sparse matrix-vector product, where the sum runs over overlapping pairs $(a,\mu)$ between the pair-product and derivative spaces once $Z$ is flattened into a vector indexed by $(a,\mu)$. Since $P^{(d)}$, $\Delta^{(d)}$, and $T_{\mathrm{mix}}^{(d)}$ only couple basis functions with overlapping support, each contracted product in \cref{eq:SequentialPushforward} is cheap, and the intermediate matrices $Y$ and $Z$ never reach the size $\left(n^{(d)}\right)^3$.\\

To combine the resulting coefficients $g_r^{(d)} \in \mathbb{R}^{\tilde n^{(d)}}$ with the two solution-space basis functions in \cref{eq:MassMatrixLowrank_1}, we define, for each $d=1,2,3$, the univariate integration matrix
\begin{equation}
        \label{eq:MassMatrixIntegration}
        \left(I_{\omega}^{(d)}\right)_{\left(\hat{i}^{(d)},\hat{j}^{(d)}\right),\alpha} = \int_0^1 \hat{\beta}^{(d)}_{\hat{i}^{(d)}}\left(\hat{x}^{(d)}\right) \hat{\beta}^{(d)}_{\hat{j}^{(d)}}\left(\hat{x}^{(d)}\right) \, \tilde{\beta}^{(d)}_{\alpha}\left(\hat{x}^{(d)}\right) \,\mathrm{d}\hat{x}^{(d)},
\end{equation}
which is of size $\left(\hat{n}^{(d)}\right)^2 \times \tilde{n}^{(d)}$ and maps coefficients in $\mathbb{S}_{3p^{(d)}-1}\left(\tilde{\Xi}^{(d)}\right)$ to flattened univariate mass matrices. Applying \cref{eq:SequentialPushforward} to every $\tilde{C}^{(d)}_{\Sigma,r}$ in \cref{eq:WeightFunctionFrobeniusProductCoefficientsTT}, followed by $I_\omega^{(d)}$, the mass tensor is assembled as
{\small
\begin{equation}
        \label{eq:LowRankMassMatrixNew}
        \mathbf{M} = \sum_{r=1}^{R_{\Sigma}} \bigotimes_{d=1}^{3} \operatorname{vec}\left(M_r^{(d)}\right) \in \mathbb{R}^{\left(\left(\hat{n}^{(1)}\right)^2,\left(\hat{n}^{(2)}\right)^2,\left(\hat{n}^{(3)}\right)^2\right)}, \quad \operatorname{vec}\left(M_r^{(d)}\right) = I_{\omega}^{(d)}\, g_r^{(d)},
\end{equation}}
with $M_r^{(d)} \in \mathbb{R}^{\hat{n}^{(d)}\times \hat{n}^{(d)}}$. Reshaping $\mathbf{M}$ finally gives\\ 
$\mathbf{M} \in \mathbb{R}^{\left(\hat{n}^{(1)},\hat{n}^{(2)},\hat{n}^{(3)}\right) \times \left(\hat{n}^{(1)},\hat{n}^{(2)},\hat{n}^{(3)}\right)}$.

We deliberately do not use an approximation symbol $\approx$ in \cref{eq:LowRankMassMatrixNew}. For a tensor-product B-spline geometry, the construction is exact provided that the tensor $\tilde{\mathbf{C}}_{\Sigma}$ is represented without tensor truncation and the arising univariate spline integrals are evaluated exactly. In this case, the weight function is represented exactly in the reduced tensor-product spline space $\mathbb{S}_{3p^{(1)}-1}\left(\tilde{\Xi}^{(1)}\right) \otimes \mathbb{S}_{3p^{(2)}-1}\left(\tilde{\Xi}^{(2)}\right) \otimes \mathbb{S}_{3p^{(3)}-1}\left(\tilde{\Xi}^{(3)}\right)$. This means that the assembled low-rank mass tensor coincides with the standard Galerkin mass tensor up to floating-point arithmetic. Approximation errors are introduced only by tensor truncation or by inexact quadrature. Therefore, if no truncation is applied, the method is an exact reformulation of the standard mass matrix assembly.

\subsection{Assembling the stiffness tensor}
\label{subsec:StiffTensor}
Next, we present a projection-based approach for assembling the stiffness tensor in low-rank format. The stiffness tensor $\mathbf{K}$, defined in \cref{eq:MassAndStiffness}, contains the matrix-valued weight function $Q$, see \cref{eq:WeightFunctions}. Although the derivatives of the tensor-product basis functions in \cref{eq:MassAndStiffness} still possess a separable structure, this alone is not sufficient to obtain a separated representation of $\mathbf{K}$, since the entries $q_{kl}$ of $Q$ are in general not separable. In contrast to the mass tensor $\mathbf{M}$, the presented construction is not exact in general, because the $q_{kl}$ are rational functions and therefore do not generally belong to a non-rational B-spline space. Nevertheless, we transfer the ideas from \cref{subsec:MassTensor} to the stiffness case. To keep the notation simple in the following, we omit the dependence on $\hat{x}$.

First, we derive a representation of the entries $q_{kl}$ in terms of Frobenius products, which allows us to exploit the underlying tensor-product structure. By \cref{eq:GeometryMapPartialDerivative}, the Jacobian of the geometry map can be written as
\begin{equation*}
        \nabla G = \left( \left\langle \mathbf{C}^{(i)}, \ \partial^{(j)}\mathbf{B} \right\rangle_F \right)_{i,j=1}^{3}.
\end{equation*}
For representing $\left(\nabla G\right)^{-1}$ in \cref{eq:WeightFunctions}, we use the adjugate matrix of $\nabla G$, i.e.
\begin{equation}
        \label{eq:AdjugateJacobianStiffness}
        \left(\nabla G\right)^{-1} = \frac{1}{\omega} \operatorname{adj}\left(\nabla G\right) = \frac{1}{\omega} \left( A_{ij} \right)_{i,j=1}^{3},
\end{equation}
where the entries $A_{ij}$ are given by
{\small
\begin{equation}
        \label{eq:AdjugateEntries}
        \begin{aligned}
                A_{ij} & = \frac{\partial G_{j+1}}{\partial \hat{x}^{(i+1)}} \frac{\partial G_{j+2}}{\partial \hat{x}^{(i+2)}} - \frac{\partial G_{j+2}}{\partial \hat{x}^{(i+1)}} \frac{\partial G_{j+1}}{\partial \hat{x}^{(i+2)}} \\
                & = \left\langle \mathbf{C}^{(j+1)}, \ \partial^{(i+1)}\mathbf{B} \right\rangle_F \left\langle \mathbf{C}^{(j+2)}, \ \partial^{(i+2)}\mathbf{B} \right\rangle_F - \left\langle \mathbf{C}^{(j+2)}, \ \partial^{(i+1)}\mathbf{B} \right\rangle_F \left\langle \mathbf{C}^{(j+1)}, \ \partial^{(i+2)}\mathbf{B} \right\rangle_F,
        \end{aligned}
\end{equation}}
where the indices are taken cyclically in $\left\{1,2,3\right\}$, i.e.~for $i=1$ and $j=2$ we have
\begin{equation*}
        A_{12} = \frac{\partial G_3}{\partial \hat{x}^{(2)}} \frac{\partial G_1}{\partial \hat{x}^{(3)}} - \frac{\partial G_1}{\partial \hat{x}^{(2)}} \frac{\partial G_3}{\partial \hat{x}^{(3)}} .
\end{equation*}
With this notation, the matrix-valued weight function can be written as
\begin{equation*}
        Q = \omega   \left( \nabla G  \right)^{-1}  \left( \nabla G  \right)^{-\top} = \rho \operatorname{adj}\left(\nabla G \right) \operatorname{adj}\left(\nabla G\right)^{\top}
\end{equation*}
and each entry has the form
\begin{equation}
        \label{eq:QEntriesNumeratorDenominator}
        q_{kl} \left( \hat{x} \right) = N_{kl} \left( \hat{x} \right) \rho \left( \hat{x} \right), \qquad k,l=1,2,3, \qquad \rho \left( \hat{x} \right) = \frac{1}{\omega \left( \hat{x} \right)},
\end{equation}
where the numerator functions are
\begin{equation}
        \label{eq:QNumeratorsStiffness}
        \begin{aligned}
                N_{11} &= A_{11}^2 + A_{12}^2 + A_{13}^2, \\
                N_{12} &= A_{11}A_{21} + A_{12}A_{22} + A_{13}A_{23}, \\
                N_{13} &= A_{11}A_{31} + A_{12}A_{32} + A_{13}A_{33}, \\
                N_{22} &= A_{21}^2 + A_{22}^2 + A_{23}^2, \\
                N_{23} &= A_{21}A_{31} + A_{22}A_{32} + A_{23}A_{33}, \\
                N_{33} &= A_{31}^2 + A_{32}^2 + A_{33}^2 ,
        \end{aligned}
\end{equation}
since $Q$ is symmetric, only these six numerator functions have to be constructed.

Considering \cref{eq:QEntriesNumeratorDenominator}, the situation is as follows. Since the entries $A_{ij}$ of the adjugate matrix $\operatorname{adj}\left(\nabla G\right)$ are sums of products of Frobenius products, as in \cref{eq:AdjugateEntries}, and since the numerator functions $N_{kl}$ are sums of products of these entries, see \cref{eq:QNumeratorsStiffness}, the numerator functions $N_{kl}$ can also be represented using Frobenius products. This representation will be derived next. As shown in \cref{subsec:MassTensor}, the denominator $\omega$ in \cref{eq:QEntriesNumeratorDenominator} is a spline function and can be represented exactly in a reduced tensor-product spline space by transferring the coefficients $\mathbf{C}^{(a)}$ to this space. However, the reciprocal $\rho$ is in general not an element of a B-spline space and therefore cannot be represented exactly by non-rational B-splines. To overcome this issue, we approximate $\rho$ by projection onto a sufficiently rich spline space. Then, by combining this approximation with the numerator functions $N_{kl}$, the entries $q_{kl}$ of $Q$ can be represented in a separated form.

To represent the numerator functions $N_{kl}$ by Frobenius products we apply similar techniques as in \cref{subsec:MassTensor}, namely the bilinearity of the Frobenius product and suitable permutations of modes. We first rewrite the entries $A_{ij}$ of the adjugate matrix \cref{eq:AdjugateJacobianStiffness} in this form. For the first row of $\operatorname{adj}\left(\nabla G\right)$, we obtain
\begin{equation*}
        A_{1j} = \left\langle \mathbf{A}_{1j}, \ \partial^{(2)}\mathbf{B} \otimes \partial^{(3)}\mathbf{B} \right\rangle_F, \qquad j = 1,2,3,
\end{equation*}
where
{\small
\begin{gather*}
        \mathbf{A}_{11} = \mathbf{C}^{(2)} \otimes \mathbf{C}^{(3)} - \mathbf{C}^{(3)} \otimes \mathbf{C}^{(2)}, \quad
        \mathbf{A}_{12} = \mathbf{C}^{(3)} \otimes \mathbf{C}^{(1)} - \mathbf{C}^{(1)} \otimes \mathbf{C}^{(3)}, \\
        \mathbf{A}_{13} = \mathbf{C}^{(1)} \otimes \mathbf{C}^{(2)} - \mathbf{C}^{(2)} \otimes \mathbf{C}^{(1)},
\end{gather*}}
are tensors of order $6$. For the entry $A_{12}$ of $\operatorname{adj}\left(\nabla G\right)$ then holds
{\small
\begin{align*}
        A_{12} & = \frac{\partial G_3}{\partial \hat{x}^{(2)}} \frac{\partial G_1}{\partial \hat{x}^{(3)}} - \frac{\partial G_1}{\partial \hat{x}^{(2)}} \frac{\partial G_3}{\partial \hat{x}^{(3)}} = \left\langle \mathbf{C}^{(3)} \otimes \mathbf{C}^{(1)} - \mathbf{C}^{(1)} \otimes \mathbf{C}^{(3)}, \ \partial^{(2)}\mathbf{B} \otimes \partial^{(3)}\mathbf{B} \right\rangle_F \\
        & = \left\langle \mathbf{A}_{12}, \ \partial^{(2)}\mathbf{B} \otimes \partial^{(3)}\mathbf{B} \right\rangle_F.
\end{align*}}

Analogously, the second and third rows can be written as
\begin{equation*}
        A_{2j} = \left\langle \mathbf{A}_{2j}, \ \partial^{(1)}\mathbf{B} \otimes \partial^{(3)}\mathbf{B} \right\rangle_F, \qquad
        A_{3j} = \left\langle \mathbf{A}_{3j}, \ \partial^{(1)}\mathbf{B} \otimes \partial^{(2)}\mathbf{B} \right\rangle_F,
\end{equation*}
with
{\small
\begin{gather*}
        \mathbf{A}_{21} = \mathbf{C}^{(3)} \otimes \mathbf{C}^{(2)} - \mathbf{C}^{(2)} \otimes \mathbf{C}^{(3)}, \quad
        \mathbf{A}_{22} = \mathbf{C}^{(1)} \otimes \mathbf{C}^{(3)} - \mathbf{C}^{(3)} \otimes \mathbf{C}^{(1)}, \\
        \mathbf{A}_{23} = \mathbf{C}^{(2)} \otimes \mathbf{C}^{(1)} - \mathbf{C}^{(1)} \otimes \mathbf{C}^{(2)}, \\
        \mathbf{A}_{31} = \mathbf{C}^{(2)} \otimes \mathbf{C}^{(3)} - \mathbf{C}^{(3)} \otimes \mathbf{C}^{(2)}, \quad
        \mathbf{A}_{32} = \mathbf{C}^{(3)} \otimes \mathbf{C}^{(1)} - \mathbf{C}^{(1)} \otimes \mathbf{C}^{(3)}, \\
        \mathbf{A}_{33} = \mathbf{C}^{(1)} \otimes \mathbf{C}^{(2)} - \mathbf{C}^{(2)} \otimes \mathbf{C}^{(1)}.
\end{gather*}
}
We point out that the tensors $\mathbf{A}_{ij}$ are of order $6$.

Using these tensors, the numerator functions can be represented as one Frobenius product. For example, since
\begin{equation*}
        N_{12} = A_{11}A_{21} + A_{12}A_{22} + A_{13}A_{23},
\end{equation*}
we obtain
\begin{equation*}
        N_{12} = \left\langle \mathbf{N}_{12}, \ \left(\partial^{(2)}\mathbf{B} \otimes \partial^{(3)}\mathbf{B}\right) \otimes \left(\partial^{(1)}\mathbf{B} \otimes \partial^{(3)}\mathbf{B}\right) \right\rangle_F,
\end{equation*}
where
\begin{equation}
        \label{eq:QNumerator12CoefficientTensor}
        \mathbf{N}_{12} = \sum_{s=1}^{3} \mathbf{A}_{1s} \otimes \mathbf{A}_{2s}.
\end{equation}
In expanded form the coefficient tensor is given by
\begin{align*}
        \mathbf{N}_{12} &= \mathbf{C}^{(2)} \otimes \mathbf{C}^{(3)} \otimes \mathbf{C}^{(3)} \otimes \mathbf{C}^{(2)} - \mathbf{C}^{(2)} \otimes \mathbf{C}^{(3)} \otimes \mathbf{C}^{(2)} \otimes \mathbf{C}^{(3)} \\
        & \quad - \mathbf{C}^{(3)} \otimes \mathbf{C}^{(2)} \otimes \mathbf{C}^{(3)} \otimes \mathbf{C}^{(2)} + \mathbf{C}^{(3)} \otimes \mathbf{C}^{(2)} \otimes \mathbf{C}^{(2)} \otimes \mathbf{C}^{(3)} \\
        & \quad + \mathbf{C}^{(3)} \otimes \mathbf{C}^{(1)} \otimes \mathbf{C}^{(1)} \otimes \mathbf{C}^{(3)} - \mathbf{C}^{(3)} \otimes \mathbf{C}^{(1)} \otimes \mathbf{C}^{(3)} \otimes \mathbf{C}^{(1)} \\
        & \quad - \mathbf{C}^{(1)} \otimes \mathbf{C}^{(3)} \otimes \mathbf{C}^{(1)} \otimes \mathbf{C}^{(3)} + \mathbf{C}^{(1)} \otimes \mathbf{C}^{(3)} \otimes \mathbf{C}^{(3)} \otimes \mathbf{C}^{(1)} \\
        & \quad + \mathbf{C}^{(1)} \otimes \mathbf{C}^{(2)} \otimes \mathbf{C}^{(2)} \otimes \mathbf{C}^{(1)} - \mathbf{C}^{(1)} \otimes \mathbf{C}^{(2)} \otimes \mathbf{C}^{(1)} \otimes \mathbf{C}^{(2)} \\
        & \quad - \mathbf{C}^{(2)} \otimes \mathbf{C}^{(1)} \otimes \mathbf{C}^{(2)} \otimes \mathbf{C}^{(1)} + \mathbf{C}^{(2)} \otimes \mathbf{C}^{(1)} \otimes \mathbf{C}^{(1)} \otimes \mathbf{C}^{(2)} .
\end{align*}
As a consequence, $\mathbf{N}_{12}$ is a tensor of order $12$, since each coordinate control-point tensor $\mathbf{C}^{(a)}$ is of order $3$. 

More generally, the numerator functions satisfy
\begin{equation*}
        N_{kl} = \sum_{s=1}^{3} A_{ks} A_{ls}, \quad k,l = 1, 2, 3.
\end{equation*}
Therefore, if
\begin{equation*}
        A_{ks} = \left\langle \mathbf{A}_{ks}, \mathbf{B}_{A,k} \right\rangle_F, \quad A_{ls} = \left\langle \mathbf{A}_{ls}, \mathbf{B}_{A,l} \right\rangle_F,
\end{equation*}
with 
\begin{equation}
        \label{eq:BSplinesPartialDerivativeTensorProduct}
        \mathbf{B}_{A,1} = \partial^{(2)}\mathbf{B} \otimes \partial^{(3)}\mathbf{B}, \quad \mathbf{B}_{A,2} = \partial^{(1)}\mathbf{B} \otimes \partial^{(3)}\mathbf{B}, \quad \mathbf{B}_{A,3} = \partial^{(1)}\mathbf{B} \otimes \partial^{(2)}\mathbf{B},
\end{equation}
then
\begin{equation}
        \label{eq:QNumeratorCoefficientGeneral}
        N_{kl} = \left\langle \mathbf{N}_{kl}, \mathbf{B}_{A,k} \otimes \mathbf{B}_{A,l} \right\rangle_F, \quad \mathbf{N}_{kl} = \sum_{s=1}^{3} \mathbf{A}_{ks} \otimes \mathbf{A}_{ls}.
\end{equation}

The representation in \cref{eq:QNumeratorCoefficientGeneral} is not yet suitable for the low-rank assembly. As for the determinant in \cref{subsec:MassTensor}, the coefficient tensors $\mathbf{N}_{kl}$ are of high order and their modes are not yet ordered according to the parametric directions.

We therefore proceed as for \cref{eq:WeightFunctionFrobeniusProductCoefficientsTT} and group the modes of each $\mathbf{N}_{kl}$ according to the three parametric directions. 

Since every summand of $\mathbf{N}_{kl}$ is a tensor product of four coordinate control-point tensors, each direction contributes a fourfold product of univariate B-splines. Following the notation of \cref{eq:SplineProduct3p-1}, only three such products occur,
{\small
\begin{equation}
        \label{eq:4pProduct}
        \begin{aligned}
                \mathrm{BBBB} \colon \quad & \left(  b^{(d)} \otimes b^{(d)} \otimes b^{(d)} \otimes b^{(d)}  \right)_{\left( i,j,k,l \right)}  = \beta_i^{\left( d \right)}  \beta_j^{\left( d \right)}  \beta_k^{\left( d \right)} \beta_l^{\left( d \right)}, \\
                \mathrm{DBBB} \colon \quad & \left(  \frac{\partial}{\partial \hat{x}^{(d)}} b^{(d)} \otimes b^{(d)} \otimes b^{(d)} \otimes b^{(d)}  \right)_{\left( i,j,k,l \right)}  = \left( \beta_i^{\left( d \right)} \right)' \beta_j^{\left( d \right)}  \beta_k^{\left( d \right)} \beta_l^{\left( d \right)}, \\
                \mathrm{DDBB} \colon \quad & \left(  \frac{\partial}{\partial \hat{x}^{(d)}} b^{(d)} \otimes \frac{\partial}{\partial \hat{x}^{(d)}} b^{(d)} \otimes b^{(d)} \otimes b^{(d)}  \right)_{\left( i,j,k,l \right)}  = \left( \beta_i^{\left( d \right)} \right)' \left( \beta_j^{\left( d \right)} \right)' \beta_k^{\left( d \right)} \beta_l^{\left( d \right)},
        \end{aligned}
\end{equation}}
depending on whether none, one, or two of the four factors are differentiated. Here we use as notation, that no occurrence gives the product type $\mathrm{BBBB}$, one occurrence gives the product type $\mathrm{DBBB}$, and two occurrences give the product type $\mathrm{DDBB}$.

Which of the three types appears in a given direction can be read off from \cref{eq:BSplinesPartialDerivativeTensorProduct}, since the entries in row $k$ of $\operatorname{adj}\left(\nabla G\right)$ carry derivatives with respect to the two directions different from $k$. It is therefore sufficient to count how often the derivative in each direction occurs in $\mathbf{B}_{A,k} \otimes \mathbf{B}_{A,l}$. For example, $N_{12}$ combines the first and second rows of the adjugate matrix. Thus, direction $3$ occurs twice, while directions $1$ and $2$ occur once, which leads to the product type $\mathrm{DBBB} \otimes \mathrm{DBBB} \otimes \mathrm{DDBB}$. As a result, the six numerator functions are assigned the following product types
\begin{equation}
        \label{eq:NumeratorProductTypeAssignment}
        \begin{array}{c|c}
        \text{numerator} & \text{product-type in dimension } (1,2,3) \\ \hline
        N_{11} & \mathrm{BBBB} \otimes \mathrm{DDBB} \otimes \mathrm{DDBB} \\
        N_{12} & \mathrm{DBBB} \otimes \mathrm{DBBB} \otimes \mathrm{DDBB} \\
        N_{13} & \mathrm{DBBB} \otimes \mathrm{DDBB} \otimes \mathrm{DBBB} \\
        N_{22} & \mathrm{DDBB} \otimes \mathrm{BBBB} \otimes \mathrm{DDBB} \\
        N_{23} & \mathrm{DDBB} \otimes \mathrm{DBBB} \otimes \mathrm{DBBB} \\
        N_{33} & \mathrm{DDBB} \otimes \mathrm{DDBB} \otimes \mathrm{BBBB}
        \end{array}
\end{equation}
The three diagonal numerators combine one product of type $\mathrm{BBBB}$ with two of type $\mathrm{DDBB}$, while the off-diagonal ones combine one of type $\mathrm{DDBB}$ with two of type $\mathrm{DBBB}$. We point out that only these three univariate product types occur at all. As a result, the coefficient transfer operators have to be constructed only once per type and parametric direction, and can then be reused for all six numerator functions.

Next, as in \cref{subsec:MassTensor}, we introduce suitable univariate reduced spline spaces that contain the products in \cref{eq:4pProduct}. If no derivative occurs, as in the first product in \cref{eq:4pProduct}, the product has degree $4p^{(d)}$. If one derivative occurs, as in the second product in \cref{eq:4pProduct}, the product has degree $4p^{(d)}-1$. If two derivatives occur, as in the third product in \cref{eq:4pProduct}, the product has degree $4p^{(d)}-2$. We denote the corresponding spaces by
\begin{equation}
        \label{eq:NumeratorSplineSpace}
        \mathbb{S}_{4p^{(d)}}\left(\Xi_{\mathrm{BBBB}}^{(d)} \right), \qquad \mathbb{S}_{4p^{(d)}-1}\left(\Xi_{\mathrm{DBBB}}^{(d)} \right), \qquad \mathbb{S}_{4p^{(d)}-2}\left(\Xi_{\mathrm{DDBB}}^{(d)} \right),
\end{equation}
where the indices $\mathrm{BBBB}$, $\mathrm{DBBB}$, and $\mathrm{DDBB}$ indicate the structure of the corresponding univariate product.

The knot vectors are constructed by the same regularity arguments as in \cref{subsec:MassTensor}. Let $\zeta^{(d)}_0 < \zeta^{(d)}_1 < \cdots < \zeta^{(d)}_{m^{(d)}}$ be the breaking points of $\Xi^{(d)}$, and let $\mu^{(d)}_\ell$ denote their multiplicities in the original geometry knot vector. For the products in \cref{eq:4pProduct}, the endpoint multiplicities are chosen as the respective degree of the reduced spline space plus $1$, while the interior multiplicities are determined by preserving the minimal regularity of the factors. We obtain the corresponding knot vectors
\begin{gather*}
        \Xi_{\mathrm{BBBB}}^{(d)} = \big\{\underbrace{\zeta^{(d)}_0,\ldots,\zeta^{(d)}_0}_{4p^{(d)}+1}, \ldots, \underbrace{\zeta^{(d)}_\ell,\ldots,\zeta^{(d)}_\ell}_{3p^{(d)}+\mu^{(d)}_\ell}, \ldots, \underbrace{\zeta^{(d)}_{m^{(d)}},\ldots,\zeta^{(d)}_{m^{(d)}}}_{4p^{(d)}+1} \big\}, \\
        \Xi_{\mathrm{DBBB}}^{(d)} = \big\{\underbrace{\zeta^{(d)}_0,\ldots,\zeta^{(d)}_0}_{4p^{(d)}}, \ldots, \underbrace{\zeta^{(d)}_\ell,\ldots,\zeta^{(d)}_\ell}_{3p^{(d)}+\mu^{(d)}_\ell}, \ldots, \underbrace{\zeta^{(d)}_{m^{(d)}},\ldots,\zeta^{(d)}_{m^{(d)}}}_{4p^{(d)}} \big\}, \\
        \Xi_{\mathrm{DDBB}}^{(d)} = \big\{\underbrace{\zeta^{(d)}_0,\ldots,\zeta^{(d)}_0}_{4p^{(d)}-1}, \ldots, \underbrace{\zeta^{(d)}_\ell,\ldots,\zeta^{(d)}_\ell}_{3p^{(d)}+\mu^{(d)}_\ell-1}, \ldots, \underbrace{\zeta^{(d)}_{m^{(d)}},\ldots,\zeta^{(d)}_{m^{(d)}}}_{4p^{(d)}-1} \big\}.
\end{gather*}
With this construction, the reduced spline spaces $\mathbb{S}_{4p^{(d)}}\left(\Xi_{\mathrm{BBBB}}^{(d)} \right)$, $\mathbb{S}_{4p^{(d)}-1}\left(\Xi_{\mathrm{DBBB}}^{(d)} \right)$ and $\mathbb{S}_{4p^{(d)}-2}\left(\Xi_{\mathrm{DDBB}}^{(d)} \right)$ contain the corresponding products \cref{eq:4pProduct}.\\

We can now transfer the coefficients of the numerator functions to the reduced spline spaces \cref{eq:NumeratorSplineSpace}. For each parametric direction $d=1,2,3$ and each product type, we introduce the univariate coefficient transfer operator
\begin{equation}
        \label{eq:FourfoldPushforward}
        F_{\mathrm{BBBB}}^{(d)} \in \mathbb{R}^{n_{\mathrm{BBBB}}^{(d)} \times \left(n^{(d)}\right)^4}, \quad F_{\mathrm{DBBB}}^{(d)} \in \mathbb{R}^{n_{\mathrm{DBBB}}^{(d)} \times \left(n^{(d)}\right)^4}, \quad F_{\mathrm{DDBB}}^{(d)} \in \mathbb{R}^{n_{\mathrm{DDBB}}^{(d)} \times \left(n^{(d)}\right)^4},
\end{equation}
which maps coefficients with respect to the fourfold products \cref{eq:4pProduct} to coefficients with respect to the basis of the corresponding reduced space. Denoting the basis of $\mathbb{S}_{4p^{(d)}}\left(\Xi_{\mathrm{BBBB}}^{(d)}\right)$ by $\left\{\beta_{\mathrm{BBBB},\alpha}^{(d)}\right\}_{\alpha=1}^{n_{\mathrm{BBBB}}^{(d)}}$, the operator is defined by
\begin{equation*}
        \beta_i^{(d)} \beta_j^{(d)} \beta_k^{(d)} \beta_l^{(d)} = \sum_{\alpha=1}^{n_{\mathrm{BBBB}}^{(d)}} \left(F_{\mathrm{BBBB}}^{(d)}\right)_{\alpha,(i,j,k,l)} \beta_{\mathrm{BBBB},\alpha}^{(d)}
\end{equation*}
and for the two remaining spaces in \cref{eq:NumeratorSplineSpace} analogously.

The operators in \cref{eq:FourfoldPushforward} are constructed from the same ingredients as in \cref{subsec:MassTensor}, namely the derivative matrix $\Delta^{(d)}$ and the $L^2$-projection of products of two univariate splines onto suitable containing spline spaces $P^{(d)}$, see \cref{eq:CoefficientTransferMatrices}. Each operator is obtained in two stages. Consider, for example, the case $\mathrm{DDBB}$. In the first stage, the two differentiated and the two undifferentiated factors are combined pairwise. The undifferentiated factors are mapped by $P^{(d)}$ into $\mathbb{S}_{2p^{(d)}}\left(\Xi_2^{(d)}\right)$, whereas the differentiated factors are first mapped by $\Delta^{(d)}$ into $\mathbb{S}_{p^{(d)}-1}\left(\Xi_\partial^{(d)}\right)$ and subsequently combined by $P^{(d)}$ in $\mathbb{S}_{2p^{(d)}-2}\left(\Xi_{2p^{(d)}-2}^{(d)}\right)$. Thus, each derivative reduces the degree of the corresponding product space by one. In the second stage, the two intermediate factors are multiplied and mapped into the target space $\mathbb{S}_{4p^{(d)}-2}\left(\Xi_{\mathrm{DDBB}}^{(d)}\right)$. This gives the intermediate spaces
\begin{equation}
        \label{eq:FourfoldStages}
        \begin{array}{c|c|c|c}
                \text{product type} & \text{first pair} & \text{second pair} & \text{target space} \\ \hline
                \mathrm{BBBB}   & \mathbb{S}_{2p^{(d)}}   & \mathbb{S}_{2p^{(d)}} & \mathbb{S}_{4p^{(d)}}\left(\Xi_{\mathrm{BBBB}}^{(d)}\right) \\
                \mathrm{DBBB}   & \mathbb{S}_{2p^{(d)}-1} & \mathbb{S}_{2p^{(d)}} & \mathbb{S}_{4p^{(d)}-1}\left(\Xi_{\mathrm{DBBB}}^{(d)}\right) \\
                \mathrm{DDBB}   & \mathbb{S}_{2p^{(d)}-2} & \mathbb{S}_{2p^{(d)}} & \mathbb{S}_{4p^{(d)}-2}\left(\Xi_{\mathrm{DDBB}}^{(d)}\right)
        \end{array}
\end{equation}
whose knot vectors are constructed by the same regularity arguments as in \cref{subsec:MassTensor}. As in \cref{eq:CoefficientTransferMatrices}, only overlapping index pairs are considered, since the corresponding products vanish otherwise. Every stage is exact up to quadrature and floating-point errors, since each target space contains by construction the products it receives.

We point out that the operators \cref{eq:FourfoldPushforward} serve as notation and are never assembled. Each of them has $\left(n^{(d)}\right)^4$ columns, one for every index tuple $(i,j,k,l)$, so that forming them would reintroduce exactly the cost the reduced spaces are meant to avoid. As for the mass tensor in \cref{eq:SequentialPushforward}, we instead apply the two stages of \cref{eq:FourfoldStages} successively and directly to the four TT cores belonging to direction $d$. This means that the mode of size $\left(n^{(d)}\right)^4$ never occurs and that the intermediate results remain of the size of the respective pair-product space. Whenever an operator in \cref{eq:FourfoldPushforward} is applied modewise in the following, this staged evaluation is understood.\\

We illustrate the modewise application of these operators for the numerator function $N_{12}$, whose coefficient tensor is given in \cref{eq:QNumerator12CoefficientTensor}. By \cref{eq:NumeratorProductTypeAssignment}, its product type is $\mathrm{DBBB} \otimes \mathrm{DBBB} \otimes \mathrm{DDBB}$. After grouping the modes of $\mathbf{N}_{12}$ according to the three parametric directions, we obtain, with a slight abuse of notation, a tensor $\mathbf{N}_{12} \in \mathbb{R}^{\left(\left(n^{(1)}\right)^4,\left(n^{(2)}\right)^4,\left(n^{(3)}\right)^4\right)}$, and the coefficients of $N_{12}$ with respect to the reduced tensor-product numerator space are obtained by
\begin{equation*}
        \tilde{\mathbf{N}}_{12} = \mathbf{N}_{12} \times_1 F_{\mathrm{DBBB}}^{(1)} \times_2 F_{\mathrm{DBBB}}^{(2)} \times_3 F_{\mathrm{DDBB}}^{(3)} \in \mathbb{R}^{\left(n_{\mathrm{DBBB}}^{(1)},n_{\mathrm{DBBB}}^{(2)},n_{\mathrm{DDBB}}^{(3)}\right)}.
\end{equation*}
This tensor contains the coefficients of $N_{12}$ in the reduced space $\mathbb{S}_{4p^{(1)}-1}\left(\Xi_{\mathrm{DBBB}}^{(1)}\right) \otimes \mathbb{S}_{4p^{(2)}-1}\left(\Xi_{\mathrm{DBBB}}^{(2)}\right) \otimes \mathbb{S}_{4p^{(3)}-2}\left(\Xi_{\mathrm{DDBB}}^{(3)}\right)$. The other numerator functions are treated analogously, where the three operators are chosen according to the product-type assignment in \cref{eq:NumeratorProductTypeAssignment}.

As for $\tilde{\mathbf{C}}_{\Sigma}$ in \cref{eq:WeightFunctionFrobeniusProductCoefficientsTT}, the reduced numerator coefficient tensors are constructed directly in TT format, starting from the TT decompositions of the coordinate control-point tensors $\mathbf{C}^{(a)}$. Every summand of $\mathbf{N}_{kl}$ in \cref{eq:QNumeratorCoefficientGeneral} is a tensor product of four such tensors, so that the four TT cores belonging to direction $d$ can be combined into a single core with TT ranks $R^{(d-1)}_{N_{kl}}$ and $R^{(d)}_{N_{kl}}$. Since the operations of \cref{eq:FourfoldPushforward} act on the physical mode only and leave the TT ranks unchanged, we apply them to each summand immediately and add the summands afterwards in the reduced spaces \cref{eq:NumeratorSplineSpace}, with a rounding after each addition. In this way, neither a tensor of order $12$ nor a mode of size $\left(n^{(d)}\right)^4$ is formed at any point. This yields the separated representation
\begin{equation}
        \label{eq:SeparatedNumeratorFunction}
        N_{kl} \left(\hat{x}\right) = \sum_{r=1}^{R_{N_{kl}}} \prod_{d=1}^{3} \tilde{N}^{(d)}_{kl,r} \cdot b_{kl}^{(d)} \left(\hat{x}^{(d)}\right), \qquad k,l=1,2,3,
\end{equation}
where the factors $\tilde{N}^{(d)}_{kl,r}$ are the low-rank factors of $\tilde{\mathbf{N}}_{kl}$ and $b_{kl}^{(d)}$ denotes the vector of basis functions of the corresponding reduced product spline space in direction $d$, again chosen according to \cref{eq:NumeratorProductTypeAssignment}. For $N_{12}$, these are the bases of the $\mathrm{DBBB}$ spaces for $d=1,2$ and of the $\mathrm{DDBB}$ space for $d=3$.

We point out that the representation \cref{eq:SeparatedNumeratorFunction} is exact up to the TT rounding, which is why no approximation symbol is used. This means that the numerator functions in \cref{eq:QEntriesNumeratorDenominator} are represented exactly in the corresponding reduced product spline spaces and, after a low-rank decomposition of their coefficient tensors, in separated form. \\

Next, we consider the remaining factor in \cref{eq:QEntriesNumeratorDenominator}, namely the reciprocal $\rho$. While the weight function $\omega$ in \cref{eq:WeightFunctionFrobeniusProduct} is a B-spline function, this is in general not the case for $\rho$. B-spline spaces are linear spaces and are therefore closed under addition and scalar multiplication, but not under taking reciprocals. On each knot span, $\omega$ is a polynomial, whereas $\rho$ is in general a rational function. To the best of the authors' knowledge, there is no separable linear space of functions which is well suited for approximating reciprocals of B-spline functions. We therefore approximate $\rho$ by projecting it onto a tensor-product B-spline space, which makes it separable and therefore representable in low-rank format. In contrast to the numerator functions $N_{kl}$, which are represented exactly, the reciprocal determinant $\rho$ in \cref{eq:QEntriesNumeratorDenominator} introduces an approximation error. 

In contrast to the interpolation approach in \cref{subsection:LowRankInterpolation}, the reciprocal $\rho$ is in the presented method approximated only once and then combined with all numerator functions $N_{kl}$. This means that we solve a single projection system instead of one system for each of the six independent entries $q_{kl}$ of $Q$. This is natural, since $\rho$ is the common non-polynomial factor in all $q_{kl}$, so that the approximation is concentrated on this one factor while the numerator functions have already been represented exactly, and the resulting representation of $\rho$ can be reused consistently for all entries of $Q$.

We denote the univariate projection spaces and the resulting tensor-product space by
\begin{equation}
        \label{eq:ProjectionSpace}
        \mathbb{S}_{\rho}^{(d)} = \mathbb{S}_{p_{\rho}^{(d)}}\left(\Xi_{\rho}^{(d)}\right), \qquad \mathbb{S}_{\rho} = \mathbb{S}_{\rho}^{(1)} \otimes \mathbb{S}_{\rho}^{(2)} \otimes \mathbb{S}_{\rho}^{(3)},
\end{equation}
with $n_{\rho}^{(d)}$ basis functions in direction $d$, and we do not fix a particular choice at this point. The construction below applies to any tensor-product B-spline space of this form, and the spaces used in our experiments are specified in \cref{sec:Numerics}. We point out, however, that this choice determines the accuracy of the whole approximation of $Q$. The space has to be rich enough to represent $\rho$, in particular the reduced regularity caused by the derivatives of the geometry map, and its richness should be comparable to that of the numerator spaces \cref{eq:NumeratorSplineSpace}, since $\rho$ is multiplied with the numerator functions afterwards. At the same time, enlarging $\mathbb{S}_{\rho}$ increases the size of the projection system, so that accuracy has to be balanced against computational cost.

Let $\mathbf{B}_{\rho} = b_{\rho}^{(1)} \otimes b_{\rho}^{(2)} \otimes b_{\rho}^{(3)}$ be the tensor-product of the basis functions of $\mathbb{S}_{\rho}$ with $b_{\rho}^{(d)} = \left[ \beta^{\left( d \right)}_{\rho, 1}, \ldots, \beta^{\left( d \right)}_{\rho, n_{\rho}^{(d)}} \right]^\top$. We approximate
\begin{equation}
        \label{eq:ReciprocalDeterminantApproximation}
        \rho(\hat{x}) \approx \rho_h(\hat{x}) = \left\langle \mathbf{D}, \mathbf{B}_{\rho}(\hat{x}) \right\rangle_F, \qquad \mathbf{D} \in \mathbb{R}^{\left(n_{\rho}^{(1)},n_{\rho}^{(2)},n_{\rho}^{(3)}\right)},
\end{equation}
where the coefficient tensor $\mathbf{D}$ is to be determined. A low-rank decomposition
\begin{equation*}
        \mathbf{D} = \sum_{r = 1}^{R_{\rho}} \bigotimes_{d = 1}^{3} D_{\rho,r}^{(d)}, \qquad D_{\rho,r}^{(d)} \in \mathbb{R}^{n_{\rho}^{(d)}},
\end{equation*}
gives the desired separable representation of $\rho_h$.

Instead of projecting the rational function $\rho$ directly, we use the identity $\omega \rho = 1$ and determine $\rho_h \in \mathbb{S}_{\rho}$ from
\begin{equation*}
        \int_{[0,1]^3} \left( \omega(\hat{x}) \rho_h(\hat{x}) - 1 \right) \left(\mathbf{B}_{\rho}(\hat{x})\right)_{\mathbf{i}} \,\mathrm{d}\hat{x} = 0
\end{equation*}
for all basis functions of $\mathbb{S}_{\rho}$. This can be interpreted as a projection of $\rho$ with respect to the weighted inner product 
\begin{equation*}
        \langle u,v \rangle_{\omega} = \int_{[0,1]^3} u v \, \omega \,\mathrm{d}\hat{x},
\end{equation*}
which is well defined since the geometry map is regular and orientation-preserving, see \cref{eq:OrientationPreserving}. Inserting \cref{eq:ReciprocalDeterminantApproximation}, we obtain the symmetric positive definite linear system
\begin{equation}
        \label{eq:ReciprocalProjectionLinearSystem}
        \begin{gathered}
                \mathbf{M}_{\rho} \cdot \mathbf{D} = \mathbf{b}_{\rho}, \\
                \left(\mathbf{M}_{\rho}\right)_{\mathbf{i},\mathbf{j}} = \int_{[0,1]^3} \beta_{\rho, \mathbf{i}} \, \beta_{\rho, \mathbf{j}} \, \omega \,\mathrm{d}\hat{x}, \quad \left(\mathbf{b}_{\rho}\right)_{\mathbf{i}} = \int_{[0,1]^3} \beta_{\rho, \mathbf{i}} \,\mathrm{d}\hat{x},
        \end{gathered}
\end{equation}
with multi-indices $\mathbf{i},\mathbf{j}$ corresponding to the tensor-product basis of $\mathbb{S}_{\rho}$. The projection matrix has the same structure as the mass tensor in \cref{eq:MassAndStiffness}, so that the heuristic from \cref{subsec:MassTensor} can be used for its low-rank assembly,
\begin{equation}
        \label{eq:ProjectionMassMatrix}
        \mathbf{M}_{\rho} = \sum_{r=1}^{R_{\Sigma}} \bigotimes_{d=1}^{3} M^{(d)}_{\rho, r},
\end{equation}
where the rank $R_{\Sigma}$ is inherited from \cref{eq:WeightFunctionFrobeniusProductCoefficientsTT}. The right-hand side is a rank-one tensor $\mathbf{b}_{\rho} = \int_{0}^{1} b^{(1)}_{\rho} \mathrm{d}\hat{x}^{(1)} \otimes \int_{0}^{1} b^{(2)}_{\rho} \mathrm{d}\hat{x}^{(2)} \otimes \int_{0}^{1} b^{(3)}_{\rho} \mathrm{d}\hat{x}^{(3)}$, since it does not contain any geometry-dependent coefficients and inherits the structure of $\mathbf{B}_{\rho}$.

The system \cref{eq:ReciprocalProjectionLinearSystem} is in general badly conditioned. We therefore change to an $L^2$-orthonormal basis before solving. For each direction, we assemble the univariate Gram matrix $M_{\rho}^{(d)}$ of $b_{\rho}^{(d)}$ (without any weight function in the integrand), which is symmetric positive definite, and compute its Cholesky factorization $M_{\rho}^{(d)} = L^{(d)} \left(L^{(d)}\right)^{\top}$. We point out that the entries of $M_{\rho}^{(d)}$ are simply given by the inner products of the univariate basis functions of $\mathbb{S}_{\rho}^{(d)}$ and should not be confused with the univariate factors $M^{(d)}_{\rho, r}$, which depend on the geometry through the low-rank representation of $\omega$. Setting $H^{(d)} = \left(L^{(d)}\right)^{-\top}$, we obtain
\begin{equation*}
        \left(H^{(d)}\right)^{\top} M_{\rho}^{(d)} H^{(d)} = I,
\end{equation*}
that is, $H^{(d)}$ maps to an $L^2$-orthonormal basis of $\mathbb{S}_{\rho}^{(d)}$.

Instead of solving \cref{eq:ReciprocalProjectionLinearSystem}, we define $\mathbf{H} = H^{(1)} \otimes H^{(2)} \otimes H^{(3)}$ and then solve
\begin{equation}
        \label{eq:TransformedLinearSystem}
        \left(\mathbf{H}^{\top} \cdot \mathbf{M}_{\rho} \cdot \mathbf{H} \right) \cdot \mathbf{Y} = \mathbf{H}^{\top} \cdot \mathbf{b}_{\rho}, \qquad \mathbf{D} = \mathbf{H} \cdot \mathbf{Y},
\end{equation}
where the transpose is applied to each univariate factor of $\mathbf{H}$. This transformation is compatible with the low-rank TT format. Since $\mathbf{H}$ is a Kronecker product of univariate matrices, it acts on each mode separately and can be applied directly to the TT cores of $\mathbf{M}_{\rho}$ and $\mathbf{b}_{\rho}$, which leaves the TT ranks unchanged. As each $H^{(d)}$ is square and invertible, the mode sizes $n_{\rho}^{(d)}$ are preserved, so that the solution tensor $\mathbf{Y}$ has the same format as $\mathbf{D}$.

For computing a solution of \cref{eq:TransformedLinearSystem}, we use the alternating minimal energy method, \amen~\cite{amen}, which minimizes an energy functional by updating the TT cores in an alternating way and thereby reduces the high-dimensional problem to a sequence of smaller linear systems. More precisely, we use \texttt{amen\_solve2.m} from the \textsc{TT-Toolbox}~\cite{tt-toolbox}. Alternatively, a TT version of the generalized minimal residual method, TT-GMRES~\cite{Dolgov:2013}, could be used.

We point out that the low-rank representation of $\mathbf{M}$ derived in \cref{subsec:MassTensor} is not mandatory for this construction, since a separated representation of the form \cref{eq:ProjectionMassMatrix} can also be assembled from the interpolation-based approach in \cref{subsection:LowRankInterpolation}.

After transforming back, we obtain the separated representation of the reciprocal
\begin{equation}
        \label{eq:Rho_H}
        \rho_h(\hat{x}) = \sum_{r=1}^{R_{\rho}} \prod_{d=1}^{3} D_{\rho,r}^{(d)} \cdot b_{\rho}^{(d)}\left(\hat{x}^{(d)}\right),
\end{equation}
where $b_{\rho}^{(d)}$ denotes the B-spline basis of $\mathbb{S}_{\rho}^{(d)}$. In contrast to the numerator representations derived above, this representation is in general not exact, since the projection onto $\mathbb{S}_{\rho}$ introduces an approximation error.

We now have all ingredients to represent the entries of $Q$ in separated form. We define
\begin{equation*}
        q_{kl,h} = N_{kl}\, \rho_h, \qquad k,l=1,2,3,
\end{equation*}
which is an approximation of $q_{kl} = N_{kl}\, \rho$. Inserting the separated representations \cref{eq:SeparatedNumeratorFunction} and \cref{eq:Rho_H} yields
\begin{equation}
        \label{eq:SeparatedProjectedQEntry}
        \begin{aligned}
                q_{kl,h} &= \sum_{r=1}^{R_{N_{kl}}} \sum_{s=1}^{R_{\rho}} \prod_{d=1}^{3} \left( \tilde{N}_{kl,r}^{(d)} \cdot b_{kl}^{(d)} \right) \left( D_{\rho,s}^{(d)} \cdot b_{\rho}^{(d)} \right) \\
                &= \sum_{r=1}^{R_{N_{kl}}} \sum_{s=1}^{R_{\rho}} \prod_{d=1}^{3} \left\langle \tilde{N}_{kl,r}^{(d)} \otimes D_{\rho,s}^{(d)}, \, b_{kl}^{(d)} \otimes b_{\rho}^{(d)} \right\rangle_F
        \end{aligned}
\end{equation}
so that before recompression the separation rank of $q_{kl,h}$ is at most $R_{N_{kl}} R_{\rho}$.

The representation \cref{eq:SeparatedProjectedQEntry} is separated with respect to the three parametric directions. In direction $d$, the integrand consists of two factors from the solution basis, one factor from the numerator space chosen according to \cref{eq:NumeratorProductTypeAssignment}, and one factor from the projection space $\mathbb{S}_{\rho}^{(d)}$. Together, this gives the univariate stiffness matrix
\begin{equation}
        \label{eq:UnivariateDirectQStiffnessMatrix}
        \left(K_{kl,rs}^{(d)}\right)_{i^{(d)},j^{(d)}} = \int_0^1 \delta(k,d)\hat{\beta}_{i^{(d)}}^{(d)} \; \delta(l,d)\hat{\beta}_{j^{(d)}}^{(d)} \; \left( \tilde{N}_{kl,r}^{(d)} \cdot b_{kl}^{(d)} \right) \, \left( D_{\rho,s}^{(d)} \cdot b_{\rho}^{(d)} \right) \, \mathrm{d}\hat{x}^{(d)},
\end{equation}
with the derivative operator $\delta$ from \cref{eq:DifferentialStiffnessOperator}. Summing the resulting tensor products over all separated terms gives
\begin{equation}
        \label{eq:QEntryStiffnessTensor}
        \mathbf{K}_{kl} = \sum_{r=1}^{R_{N_{kl}}} \sum_{s=1}^{R_{\rho}} K_{kl,rs}^{(1)} \otimes K_{kl,rs}^{(2)} \otimes K_{kl,rs}^{(3)},
\end{equation}
and the full stiffness tensor is obtained by summing over all entries of $Q$,
\begin{equation}
        \label{eq:LowRankStiffnessNew}
        \mathbf{K} = \sum_{k=1}^{3} \sum_{l=1}^{3} \mathbf{K}_{kl}.
\end{equation}

We point out that \cref{eq:SeparatedProjectedQEntry,eq:QEntryStiffnessTensor} describe $q_{kl,h}$ and $\mathbf{K}_{kl}$ only mathematically. In the implementation, the matrices $K_{kl,rs}^{(d)}$ are never formed term by term. Instead, the univariate TT factors of $N_{kl}$ and $\rho_h$ are contracted directly against the one-dimensional quadrature in each direction, so that neither the product coefficient tensor of $N_{kl}\rho_h$ nor the individual rank-one contributions are assembled explicitly. A further consequence of \cref{eq:UnivariateDirectQStiffnessMatrix} is that the analysis basis enters only through the two factors $\delta(k,d)\hat{\beta}_{i^{(d)}}^{(d)}$ and $\delta(l,d)\hat{\beta}_{j^{(d)}}^{(d)}$. Once the separated representations of $N_{kl}$ and $\rho_h$ are available, refining the solution space therefore only requires reassembling these univariate integrals, while the coefficient representations remain unchanged.

\section{Limitations and implementation details}
\label{sec:Limitations}

The low-rank formulas \cref{eq:LowRankMassMatrixNew,eq:LowRankStiffnessNew} hide three properties that determine the practical performance of the method, which we discuss before turning to the implementation.

The first is structural. Both constructions rely on the geometry map being an orientation-preserving tensor-product B-spline map, see \cref{eq:GeometryMapTensor}, since only then do the weight function $\omega$ and the numerator functions $N_{kl}$ separate into univariate factors. Geometries that are not of this form, in particular not tensor-product or rational parameterizations such as NURBS, are therefore outside the scope of this paper.

The second limitation concerns the coordinate control-point tensors $\mathbf{C}^{(a)} \in \mathbb{R}^{\left(n^{(1)}, n^{(2)}, n^{(3)}\right)}$, $a=1,2,3$, both through their size and through their low-rank property. The determinant tensor $\tilde{\mathbf{C}}_{\Sigma}$ in \cref{eq:WeightFunctionFrobeniusProductCoefficientsLowRank} is built from tensor products of three of these tensors and is of order nine before grouping, and the numerator tensors $\mathbf{N}_{kl}$ in \cref{eq:QNumeratorCoefficientGeneral} from tensor products of four and are therefore of order twelve. Neither is ever formed explicitly, but their TT ranks grow with those of the $\mathbf{C}^{(a)}$ and with the products and mode permutations involved, so that the method is efficient precisely when the geometry admits control-point tensors of moderate TT rank. Independently of the ranks, a large number of univariate basis functions is a bottleneck in itself. The univariate integrals and coefficient transfers scale linearly in $n^{(d)}$ for fixed degree, but a geometry with, say, $n^{(d)} = 10^4$ control points per parametric direction already leads to large univariate operators and correspondingly costly one-dimensional integration.

One possible way to address this is a low-rank multi-patch approach as in~\cite{Riemer2025}, in which the geometry is decomposed into several patches and the low-rank assembly is applied patchwise. This can reduce both the number of control points per patch and the TT ranks of the individual patch geometries, at the price of the additional coupling and interface treatment that a multi-patch formulation entails. An alternative is a targeted refinement of the spline space of the geometry map, with the aim of an equivalent representation whose control-point tensors are of lower rank. However, there is no guarantee that such a refinement reduces the ranks in practice, as this depends on the structure of the geometry.

The third limitation is the projection system \cref{eq:ReciprocalProjectionLinearSystem} for the reciprocal determinant $\rho$. Its accuracy is governed by the richness of the projection space $\mathbb{S}_{\rho}$, which has to be chosen fine enough to approximate the rational function $\rho$ well. As already discussed, improving this approximation requires knot refinement or degree elevation of $\mathbb{S}_{\rho}$, which enlarges both the projection matrix and the right-hand side and increases the degree of the univariate integrands. This can be also a challenge for low-rank solvers such as \amen. Even though the assembly and the solve are carried out in low-rank format, the univariate computations do not disappear and can become the dominant cost. A remedy would be a projection space that still has tensor-product structure but approximates $\rho$ to comparable accuracy with fewer basis functions, or the use of alternative low-rank solvers and preconditioners better suited to the weighted projection system than \amen.

Next, we describe the implementation details that influence the resulting TT ranks and the accuracy. For both assembly methods, the coordinate control-point tensors $\mathbf{C}^{(a)}$, $a = 1, 2, 3$, in \cref{eq:GeometryMapTensor} are converted into TT format using \texttt{tt\_tensor} from the \textsc{TT-Toolbox}~\cite{tt-toolbox} with tolerance \texttt{tol}. For the mass tensor assembly in \cref{subsec:MassTensor}, the determinant coefficient tensor $\tilde{\mathbf{C}}_{\Sigma}$ is then assembled directly in TT format according to \cref{eq:WeightFunctionFrobeniusProductCoefficients}, with rounding at tolerance \texttt{tol} after the tensor products, after the signed summation of the six determinant terms, and after the permutation and reshaping into the grouped tensor \cref{eq:WeightFunctionFrobeniusProductCoefficientsTT}, in order to control the intermediate TT ranks. 

For the mass tensor, the coefficient transfer in each direction is realized by the successive application of the univariate derivative matrix $\Delta^{(d)}$, obtained from the B-spline derivative formula \cref{eq:BSplineDerivative}, the univariate product transfer matrix $P^{(d)}$, and the mixed product transfer matrix $T_{\mathrm{mix}}^{(d)}$, whose coefficients are computed by one-dimensional $L^2$-projections. All univariate integrals are evaluated by Gaussian quadrature with a number of points chosen to integrate the respective polynomial integrand exactly, and entries below $10^{-14}$ are discarded in the transfer matrices to preserve sparsity. These operators, together with the univariate integration matrices $I_{\omega}^{(d)}$ in \cref{eq:LowRankMassMatrixNew}, are applied directly to the TT cores of the direction they belong to, so that the reduced coefficient tensors are never stored. After each modewise application, the tensor is rounded with tolerance \texttt{tol}. Exact integration of \cref{eq:MassMatrixIntegration} requires
\begin{equation}
        \label{eq:MassQuadrature}
        \left\lceil \frac{2\hat{p}^{(d)}+3p^{(d)}}{2} \right\rceil
\end{equation}
quadrature points in direction $d = 1, 2, 3$.

For the stiffness-tensor assembly in \cref{subsec:StiffTensor}, the matrix $\mathbf{M}_{\rho}$ is assembled in the same low-rank way as the mass tensor, and its right-hand side is a rank-one tensor of univariate integrals. Both are rounded at tolerance \texttt{tol}$/10$. When assembling \cref{eq:TransformedLinearSystem}, if $M_{\rho}^{(d)}$ becomes numerically only positive semidefinite for a strongly refined projection space and the Cholesky factorization fails, we fall back to a symmetric eigendecomposition and discard the dependent directions whose eigenvalues lie below a relative threshold. The system is then solved by \texttt{amen\_solve2.m} from the \textsc{TT-Toolbox}~\cite{tt-toolbox} with tolerance \texttt{tol}, and the resulting tensor $\mathbf{D}$ is rounded at tolerance \texttt{tol}$/10$.

The numerator coefficient tensors $\mathbf{N}_{kl}$ in \cref{eq:QNumeratorCoefficientGeneral} are assembled analogously to the determinant coefficient tensor, but are never formed explicitly. As for the mass tensor, the coefficient transfer is realized by the successive application of $\Delta^{(d)}$, $P^{(d)}$, and the mixed transfer matrices, again with exact Gaussian quadrature and entries below $10^{-14}$ discarded. The summands, the modewise applications, and the intermediate sums are rounded at tolerance \texttt{tol}$/10$. Finally, the product of $N_{kl}$ and $\rho_h$ in \cref{eq:SeparatedProjectedQEntry} is not assembled either, but contracted directly in the univariate integrals \cref{eq:UnivariateDirectQStiffnessMatrix}. Exact integration of the integrands involving $\rho_h$ requires
\begin{equation}
        \label{eq:StiffnessQuadrature}
        \left\lceil \frac{2\hat{p}^{(d)} + p_{\rho}^{(d)} + 4p^{(d)} + 1}{2} \right\rceil
\end{equation}
quadrature points in direction $d = 1, 2, 3$, where $4p^{(d)}$ bounds the numerator degrees. Each summand in \cref{eq:LowRankStiffnessNew} is rounded at tolerance \texttt{tol}$/10$, and the final stiffness tensor $\mathbf{K}$ once more at the target tolerance \texttt{tol}.

\section{Numerical experiments}
\label{sec:Numerics}

In this section, we present the results of our numerical experiments, which assess to what extent the proposed methods can compete with the full assembly implemented in \textsc{GeoPDEs} and with the low-rank interpolation approach of~\cite{BuengerDolgovStoll:2020}. We consider the mass tensor first and then the stiffness tensor.

\begin{figure}[htbp]
\centering
\begin{subfigure}[t]{0.225\textwidth}
\centering
\includegraphics[width=\linewidth]{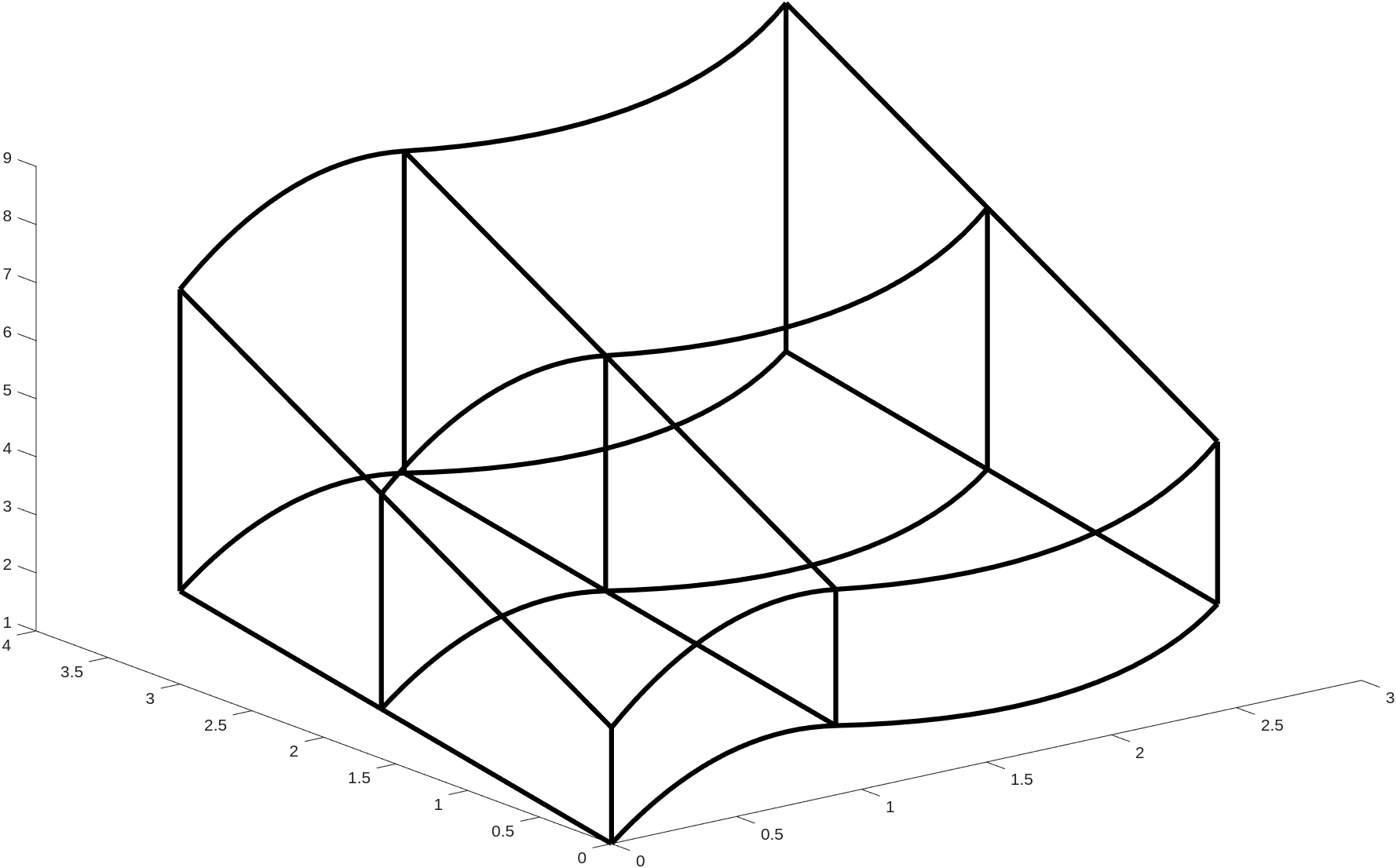}
\caption{\emph{Thick} flag}
\label{fig:1_a}
\end{subfigure}\hfill
\begin{subfigure}[t]{0.225\textwidth}
\centering
\includegraphics[width=0.2\linewidth]{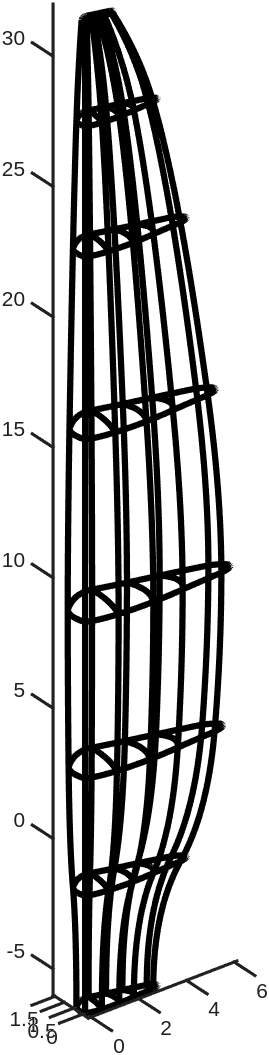}
\caption{Rotor blade}
\label{fig:1_b}
\end{subfigure}\hfill
\begin{subfigure}[t]{0.225\textwidth}
\centering
\includegraphics[width=\linewidth]{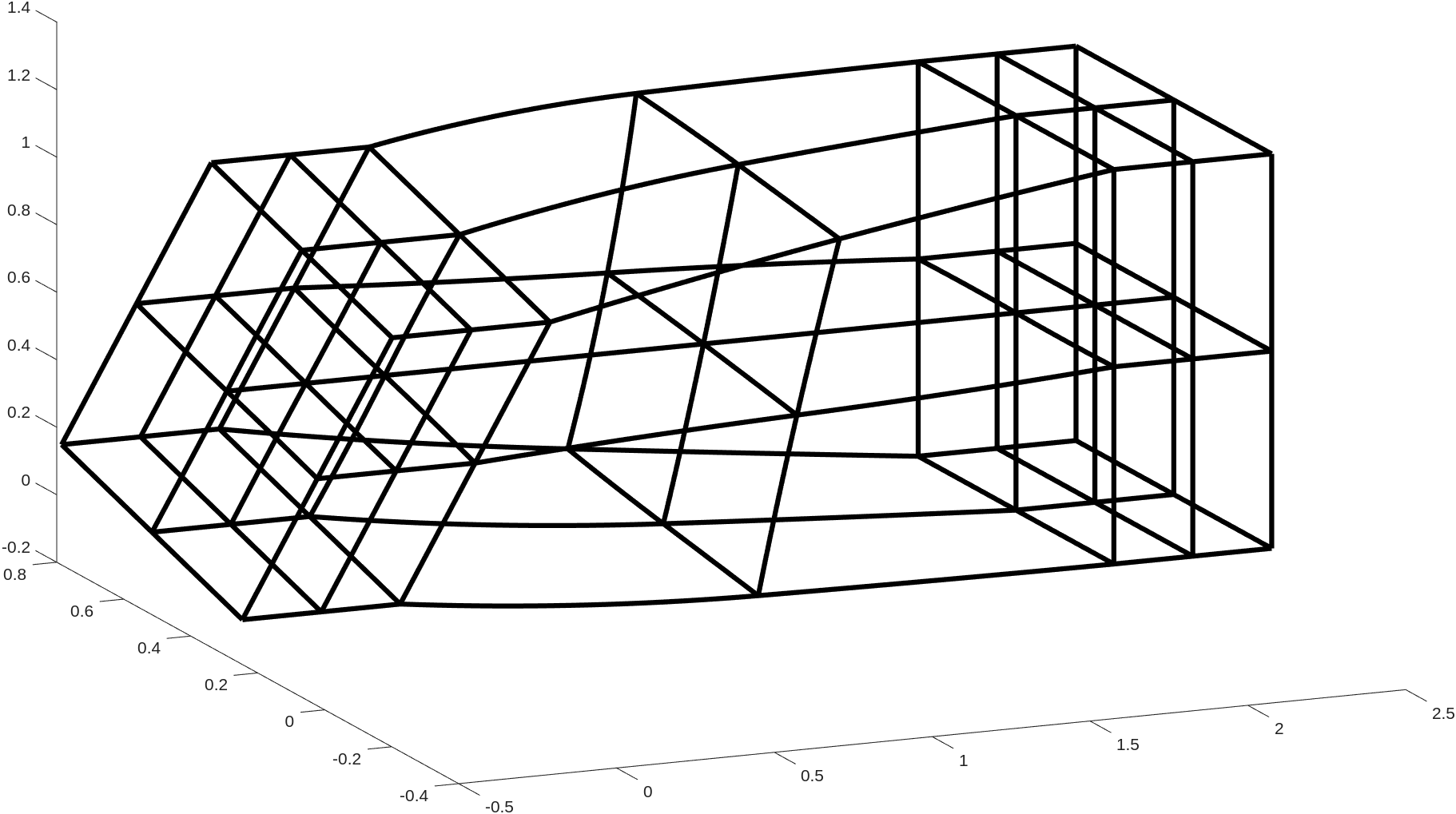}
\caption{Twisted beam}
\label{fig:1_c}
\end{subfigure}\hfill
\begin{subfigure}[t]{0.225\textwidth}
\centering
\includegraphics[width=\linewidth]{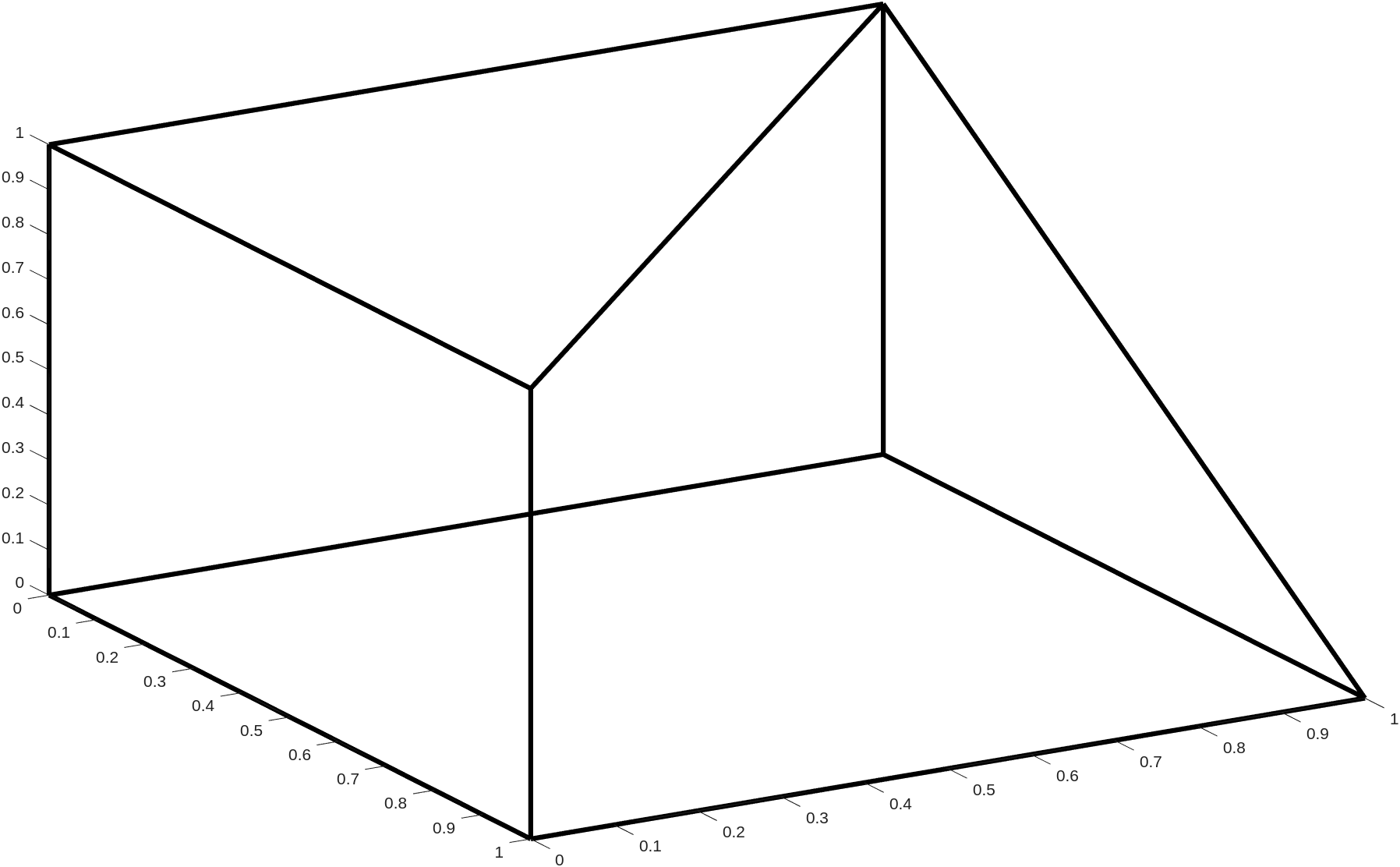}
\caption{Almost singular geometry}
\label{fig:1_d}
\end{subfigure}
\caption{Geometries used in the numerical experiments.}
\label{fig:1}
\end{figure}

The mass tensor experiments are carried out for three orientation-preserving tensor-product B-spline geometries in the sense of \cref{eq:OrientationPreserving}, which represent different levels of difficulty for the low-rank approximation. The first is a \emph{thick} flag-like volume (\cref{fig:1_a}) of degree $\mathbf{p}=(2,2,2)$ with $\mathbf{n}=(4,4,3)$ univariate basis functions. The second is a rotor blade (\cref{fig:1_b}) of degree $\mathbf{p}=(2,2,2)$ with $\mathbf{n}=(3,6,9)$. The third is a twisted beam taken from \textsc{GeoPDEs} (\cref{fig:1_c}) of degree $\mathbf{p}=(2,2,2)$ with $\mathbf{n}=(10,4,4)$. For the stiffness tensor we consider one additional geometry (\cref{fig:1_d}), for which the geometry map becomes almost singular, i.e.\ the determinant of its Jacobian becomes very small. It is of degree $\mathbf{p}=(1,1,1)$ with $\mathbf{n}=(2,2,2)$, and its construction is explained in detail at the beginning of \cref{subsec:NumericsStiffTensor}.

In each experiment, we assemble the respective matrix or tensor for solution-space degrees $\hat{p}^{(d)} = 3$ and $\hat{p}^{(d)} = 5$, $d = 1, 2, 3$, increase the number of solution basis functions by uniform refinement, and consider different tolerances for the low-rank approximation. 

In both low-rank methods, a single prescribed tolerance $\texttt{tol}$ is used which we vary it in the experiments. For the interpolation-based \cref{subsection:LowRankInterpolation} method, it serves as the stopping tolerance of \texttt{amen\_block\_solve.m} for the interpolation system \cref{eq:InterpolationSystem}, which is called with $\texttt{nswp} = 20$, $\texttt{kickrank} = 2$, $\texttt{resid\_damp} = 10$ and $\texttt{exitdir} = -1$, and as the rounding tolerance for the resulting coefficient tensor and for the partial sums of the final assembly, for the construction of the mass tensor and the stiffness tensor. For the projection-based assembly of the mass tensor, $\texttt{tol}$ is the tolerance of the final rounding of $\mathbf{M}$ as well for all intermediate roundings of the mass assembly. For the stiffness tensor $\texttt{tol}$ is the stopping tolerance of \texttt{amen\_solve2.m} for the projection system \cref{eq:TransformedLinearSystem}, called with $\texttt{nswp} = 1000$, $\texttt{kickrank} = 2$, $\texttt{resid\_damp} = 10$, $\texttt{trunc\_norm} = \texttt{'residual'}$ and $\texttt{max\_full\_size} = 100$, it is the tolerance of the final rounding of $\mathbf{K}$ and all intermediate roundings of the stiffness assembly are performed at $\texttt{tol}/10$, so that the accumulation of rounding errors does not consume the prescribed target accuracy. Only the coordinate control-point tensors $\mathbf{C}^{(a)}$ and the projection matrix $\mathbf{M}_{\rho}$ of the stiffness assembly are rounded at the fixed tolerances $10^{-13}$ and $10^{-12}$, since the accuracy of the projection operator should not be limited by the target accuracy of the assembly.

We use two quadrature rules in the final assembly step. Five Gaussian points per parametric direction as a common practical choice, and the numbers specified in \cref{eq:MassQuadrature} for the mass tensor and \cref{eq:StiffnessQuadrature} for the stiffness tensor, which integrate the polynomial integrands in \cref{eq:MassMatrixIntegration,eq:UnivariateDirectQStiffnessMatrix} exactly. The number of quadrature points used per direction is stated explicitly in each experiment. \\

For all methods we compare the memory requirements and the assembly times, and for the two low-rank methods we additionally report the relative error with respect to the fully assembled reference matrix obtained from \textsc{GeoPDEs}.

We estimate the memory required to store an assembled operator using a recursive size estimator that traverses the corresponding data structures, including sparse matrices and tensors in TT format, and accumulates the memory of all contained components. Shared references are counted only once to avoid double counting. Since low-level runtime overhead and copy-on-write effects are not accessible from user-level code, the resulting values provide an approximation of the actual memory consumption. We use this estimate to compare the storage requirements of the mass and stiffness matrices assembled with \textsc{GeoPDEs} with those of the corresponding low-rank tensors.

In addition to the operator storage, we measure the peak memory required during the assembly, since temporary quantities may require substantially more memory than the final operator. For this purpose, each configuration is executed in a separate process and we record the maximum resident set size reported by the operating system. The measured peak is taken relative to a baseline obtained after the geometry and the discrete spaces have been initialized and before the assembly is started. Running each configuration in a separate process avoids dependencies on previously allocated memory retained by the runtime and therefore provides a more reliable estimate of the transient memory required by a single assembly. The measurement is performed at the process level and consequently contains a fixed overhead common to all methods, which does not affect their relative comparison.

The assembly times are measured with \texttt{tic} and \texttt{toc}. For \textsc{GeoPDEs}, we measure only the time to assemble the reference matrix. For the interpolation approach of~\cite{BuengerDolgovStoll:2020}, the reported time includes both the determination of the interpolation coefficients and the low-rank assembly. For the proposed method, it includes the complete low-rank assembly process.

For each fixed degree, quadrature rule and refinement level, the \textsc{GeoPDEs} reference matrix $A_{\mathrm{ref}}$ is assembled once and reused for the error computation of all low-rank approximations of that configuration. If $A_{\mathrm{LR}}$ denotes the matrix represented by a low-rank tensor, we define the relative Frobenius error by
\begin{equation*}
        \varepsilon_{\mathrm{rel}} = \frac{\left\| A_{\mathrm{LR}} - A_{\mathrm{ref}} \right\|_F}{\left\| A_{\mathrm{ref}} \right\|_F},
\end{equation*}
where the low-rank tensor is converted to a sparse matrix only for the error evaluation. \\

The experiments were implemented in \textsc{MATLAB R2025b} using the \textsc{TT-Toolbox}~\cite{tt-toolbox} version 2.2.2, the \textsc{Tensor Toolbox for MATLAB}~\cite{tensor_toolbox} version 3.6, the \textsc{NURBS} package~\cite{nurbs_toolbox} version 1.4.3, the \textsc{GeoPDEs} package~\cite{geopdes3.0,geopdes_new}. 

The experiments were run on an Ubuntu Linux server (Linux kernel~5.15.0-41-generic, x86\_64) with two AMD EPYC 9534 processors ($2\times64$~cores, 2~threads per core, 256~logical CPUs total) with a maximum clock frequency of approximately $3.72\,\mathrm{GHz}$ and $1.5\,\mathrm{TiB}$ RAM.\footnote{The code is public at \url{https://github.com/TomRiem/Low_Rank_IGA_Projection.git}.}

\subsection{Mass tensor}
\label{subsec:NumericsMassTensor}

\begin{figure}[tp]
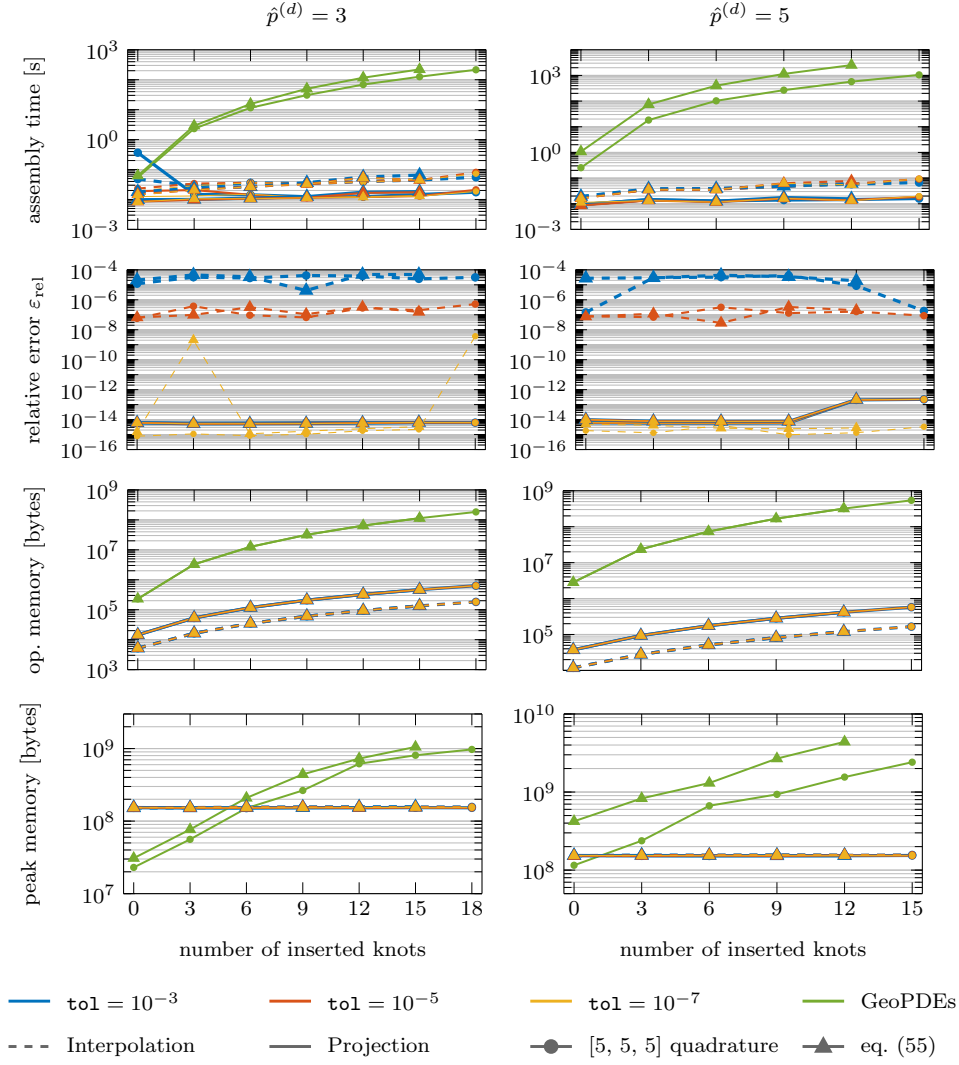

  \centering
  \setlength{\tabcolsep}{2pt}
  \begin{tabular}{cc}
    \input{flag_r1c1.tex} & \input{flag_r1c2.tex}\\
    \input{flag_r2c1.tex} & \input{flag_r2c2.tex}\\
    \input{flag_r3c1.tex} & \input{flag_r3c2.tex}\\
    \input{flag_r4c1.tex} & \input{flag_r4c2.tex}\\
  \end{tabular}\\[4pt]
  \input{flag_legend.tex}
  \caption{Results of the mass matrix/tensor assembly for the \emph{thick} flag geometry \cref{fig:1_a} for $\hat{p}^{(d)} = 3$ and $\hat{p}^{(d)} = 5$ (columns).}
  \label{fig:flag_mass}
\end{figure}

\begin{figure}[tp]
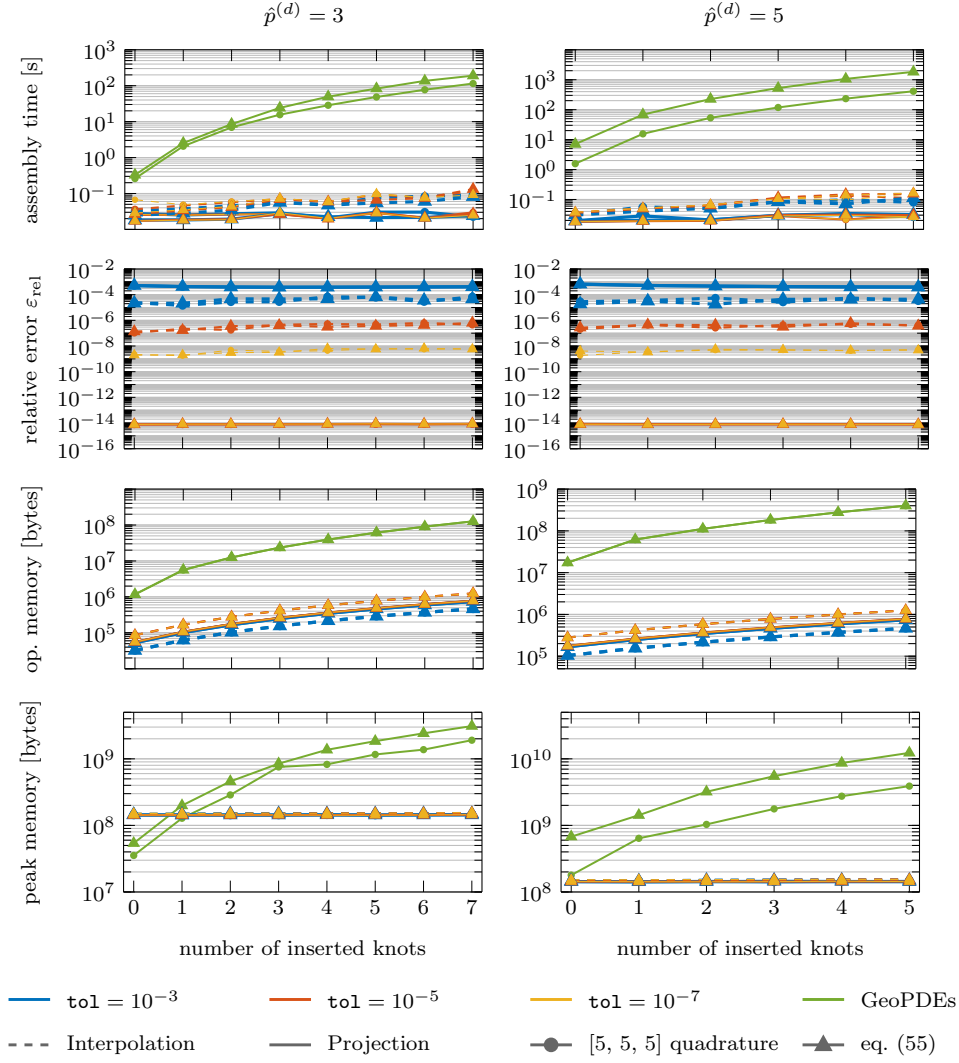

  \centering
  \setlength{\tabcolsep}{2pt}
  \begin{tabular}{cc}
    \input{rotor_r1c1.tex} & \input{rotor_r1c2.tex}\\
    \input{rotor_r2c1.tex} & \input{rotor_r2c2.tex}\\
    \input{rotor_r3c1.tex} & \input{rotor_r3c2.tex}\\
    \input{rotor_r4c1.tex} & \input{rotor_r4c2.tex}\\
  \end{tabular}\\[4pt]
  \input{rotor_legend.tex}
  \caption{Results of the mass matrix/tensor assembly for the rotor geometry \cref{fig:1_b} for $\hat{p}^{(d)} = 3$ and $\hat{p}^{(d)} = 5$ (columns).}
  \label{fig:rotor_mass}
\end{figure}

\begin{figure}[tp]
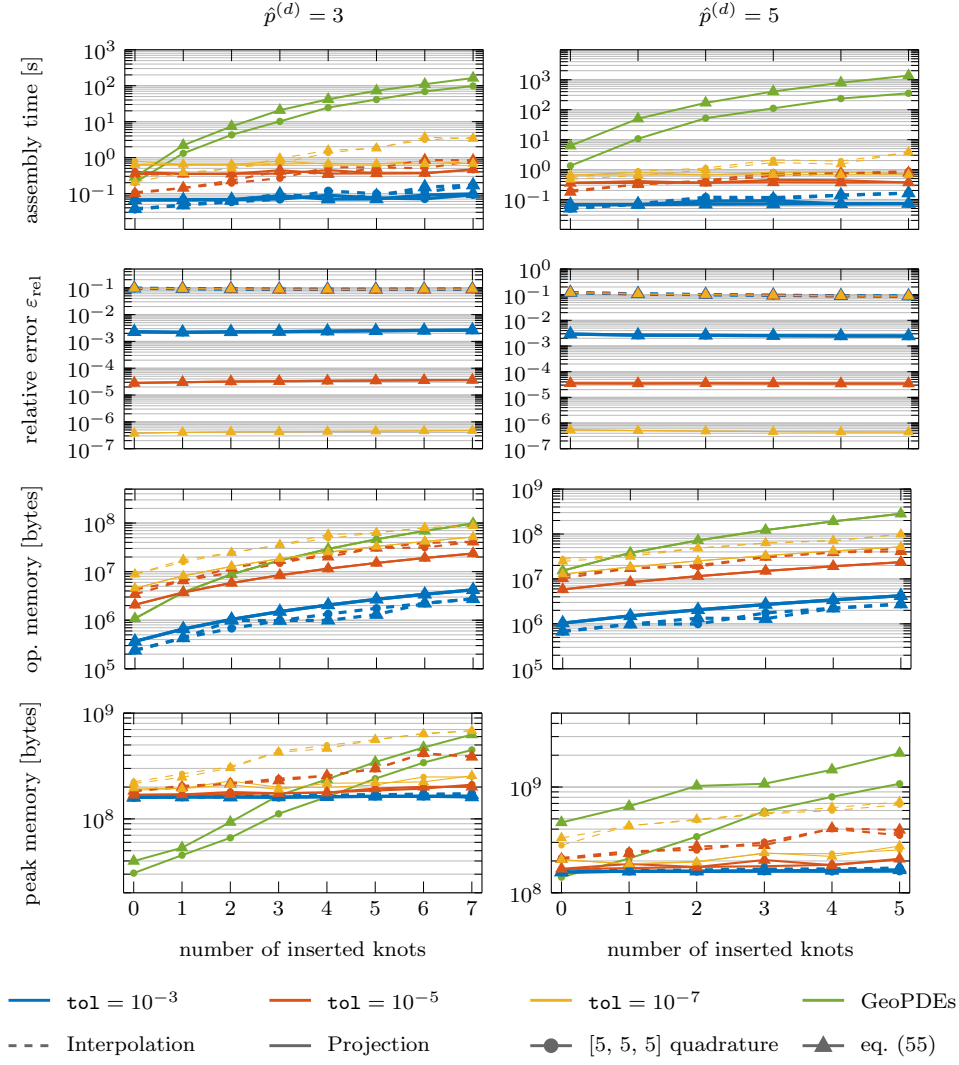

  \centering
  \setlength{\tabcolsep}{2pt}
  \begin{tabular}{cc}
    \input{beam_r1c1.tex} & \input{beam_r1c2.tex}\\
    \input{beam_r2c1.tex} & \input{beam_r2c2.tex}\\
    \input{beam_r3c1.tex} & \input{beam_r3c2.tex}\\
    \input{beam_r4c1.tex} & \input{beam_r4c2.tex}\\
  \end{tabular}\\[4pt]
  \input{beam_legend.tex}
  \caption{Results of the mass matrix/tensor assembly for the twisted beam geometry \cref{fig:1_c} for $\hat{p}^{(d)} = 3$ and $\hat{p}^{(d)} = 5$ (columns).}
  \label{fig:beam_mass}
\end{figure}

The results of the mass tensor assembly are summarized in \cref{fig:flag_mass} for the \emph{thick flag} geometry, in \cref{fig:rotor_mass} for the rotor blade, and in \cref{fig:beam_mass} for the twisted beam.\\

The results in \cref{fig:flag_mass} show that the proposed projection-based method is consistently faster than the interpolation-based assembly following~\cite{BuengerDolgovStoll:2020}, while the resulting mass tensor requires slightly more memory. However, the projection-based assembly yields errors close to machine precision for all considered tolerances, indicating that the coefficient tensor $\tilde{\mathbf{C}}_{\Sigma}$ is already approximated sufficiently well for $\varepsilon=10^{-3}$ and tighter tolerances do not provide a relevant improvement. The less regular error behavior of the interpolation-based assembly can be attributed to numerical rounding effects, while the error remains below the prescribed tolerance in all cases. For the experiments using the higher quadrature rule, we perform one refinement step less due to the increased computational time, since the available results already confirm the observed behavior. \\

For the rotor blade, interpolation-based assembly yields errors below the prescribed tolerance, the projection-based method reaches almost machine precision already for $\varepsilon=10^{-5}$. \\

The twisted beam represents the most challenging of the considered geometries. As shown in \cref{fig:beam_mass}, the interpolation-based assembly fails to approximate the mass tensor accurately, with a relative error of approximately $10^{-1}$ that is essentially independent of the prescribed tolerance. In contrast, the proposed projection-based method follows the prescribed tolerance over all refinement steps, although the errors are slightly above the respective thresholds. These small deviations can be attributed to the accumulation of TT-rounding errors. The example is also more demanding in terms of memory than the previous geometries, indicating to a more complex low-rank representation.

The failure of the interpolation-based assembly of~\cite{BuengerDolgovStoll:2020} can be explained by the reduced regularity of the twisted beam geometry of \cref{fig:1_c}. In the first parametric direction, its open knot vector is
\begin{equation*}
        \Xi^{(1)} = \left[ 0,0,0,\tfrac{1}{6},\tfrac{1}{3},\tfrac{1}{3},\tfrac{1}{2},\tfrac{2}{3},\tfrac{2}{3},\tfrac{5}{6},1,1,1 \right],
\end{equation*}
so that the interior knots $\zeta_2 = \tfrac{1}{3}$ and $\zeta_4 = \tfrac{2}{3}$ have multiplicity $\mu^{(1)}_\ell = 2 = p^{(1)}$. At such knots the geometry map is only $C^0$, and the weight function $\omega$, which involves first derivatives of the geometry, is therefore discontinuous in the first parametric direction.

By the regularity argument of \cref{subsec:MassTensor}, $\omega$ lies in this direction in the auxiliary space $\mathbb{S}_{3p^{(1)}-1}\left(\tilde{\Xi}^{(1)}\right)$ of degree $3p^{(1)}-1 = 5$, whose interior knot multiplicities are $\tilde{\mu}^{(1)}_\ell = 2p^{(1)} + \mu^{(1)}_\ell$, cf.~\cref{eq:ReducedTripleProductKnotVector}. At the two double knots this gives $\tilde{\mu}^{(1)}_\ell = 2p^{(1)} + p^{(1)} = 3p^{(1)} = 6$, i.e.\ full multiplicity for a spline of degree $3p^{(1)}-1$, so that the auxiliary space is itself $C^{-1}$ there, consistent with the discontinuity of $\omega$. As a consequence, the two double knots each give rise to a repeated Greville point, so that the entries of $\tilde{X}^{(1)}$ are no longer distinct. The univariate factor matrix $\tilde{b}^{(1)}\left(\tilde{X}^{(1)}\right)$ therefore has two pairs of identical rows, and the tensor-product interpolation system \cref{eq:InterpolationSystem} becomes singular. Here this results in a rank deficiency of two, one for each double knot, so that the interpolation coefficients cannot be determined and the interpolation-based method \cref{subsection:LowRankInterpolation} is not applicable.

For the proposed method \cref{subsec:MassTensor}, this reduced continuity is not an issue. The coefficient transfer operators $P^{(d)}$ and $T_{\mathrm{mix}}^{(d)}$ (and likewise the projection of the reciprocal determinant \cref{eq:ReciprocalProjectionLinearSystem}), are defined by $L^2$-projection onto the same, discontinuous, auxiliary spline spaces, rather than by interpolation at Greville points. Such a projection requires only the Gram matrix of the target basis, which is symmetric positive definite for any knot multiplicity, since a B-spline basis remains linearly independent even at full-multiplicity knots. The corresponding univariate integrals are evaluated by Gaussian quadrature applied separately on each knot span, on which the basis functions are polynomials and the quadrature nodes lie strictly in the interior. The quadrature therefore integrates between consecutive knot values and never across a knot, so that a discontinuity located at a repeated knot is confined to a single point and does not affect the integrals. Consequently, the projection systems remain well posed, and the discontinuity of $\omega$ at the $C^0$ knots of the geometry is no difficulty for the proposed method.

\subsection{Stiffness tensor}
\label{subsec:NumericsStiffTensor}

It remains to specify the projection space $\mathbb{S}_{\rho}$ \cref{eq:ProjectionSpace} onto which the reciprocal determinant is projected in \cref{eq:ReciprocalProjectionLinearSystem}. To make the comparison with the interpolation-based approach fair, we always use the same B-spline space for the interpolation in \cref{subsection:LowRankInterpolation} as for the projection of $\rho$. 

We use two heuristics to determine the projection space. As a default, we choose in each parametric direction a univariate spline space that contains the sixfold products
\begin{equation}
        \label{eq:SixfoldProduct}
        \beta_i^{(d)} \beta_j^{(d)} \beta_k^{(d)} \beta_l^{(d)} \left(\beta_m^{(d)}\right)' \left(\beta_n^{(d)}\right)'
\end{equation}
These products consist of four B-splines of degree $p^{(d)}$ and two derivatives of degree $p^{(d)}-1$ and therefore have degree $6p^{(d)}-2$. Following the regularity argument of \cref{subsec:MassTensor}, we construct the corresponding knot vector $\Xi_{6p^{(d)}-2}^{(d)}$ with endpoint multiplicity $6p^{(d)}-1$ and, at each interior knot $\zeta_\ell^{(d)}$ of multiplicity $\mu_\ell^{(d)}$ in $\Xi^{(d)}$, multiplicity $5p^{(d)}+\mu_\ell^{(d)}-1$. The resulting space therefore reproduces the reduced continuity $C^{p^{(d)}-\mu_\ell^{(d)}-1}$ induced by the derivative factors. We define
\begin{equation}
        \label{eq:DefaultProjectionSpace}
        \begin{gathered}
                \mathbb{S}_{\rho}^{(d)} = \mathbb{S}_{6p^{(d)}-2}\left(\Xi_{6p^{(d)}-2}^{(d)}\right), \qquad \mathbb{S}_{\rho} = \mathbb{S}_{\rho}^{(1)} \otimes \mathbb{S}_{\rho}^{(2)} \otimes \mathbb{S}_{\rho}^{(3)}, \\
                \Xi_{6p^{(d)}-2}^{(d)} = \big\{ \underbrace{\zeta^{(d)}_0,\ldots,\zeta^{(d)}_0}_{6p^{(d)}-1}, \ldots, \underbrace{\zeta^{(d)}_\ell, \ldots, \zeta^{(d)}_\ell}_{5p^{(d)}+\mu_\ell^{(d)}-1}, \ldots, \underbrace{\zeta^{(d)}_{m^{(d)}},\ldots,\zeta^{(d)}_{m^{(d)}}}_{6p^{(d)}-1} \big\}.
        \end{gathered}
\end{equation}
The sixfold product in \cref{eq:SixfoldProduct} represents the most demanding of the product structures considered here, so that the four product types in \cref{eq:4pProduct} are also contained in this space. Moreover, the chosen knot multiplicities reproduce the reduced regularity of both the determinant and the derivative-dependent numerator functions. A smaller space without derivative factors may be sufficient for particularly smooth geometries, but does not in general account for this loss of regularity. We therefore use $\mathbb{S}_{\rho}$ as a robust default projection space rather than as an optimal choice.

As a second heuristic, we use for $\mathbb{S}_{\rho}$ the tensor-product determinant space of \cref{eq:WeightFunctionSplineSpace}, but generated from the solution space $\hat{\mathbb{S}}_{\hat{p}}\left(\hat{\Xi}^{(1)}, \hat{\Xi}^{(2)}, \hat{\Xi}^{(3)}\right)$ instead of the geometry space. That is, in each parametric direction we build a B-spline space $\mathbb{S}_{\rho}^{(d)} = \mathbb{S}_{3\hat{p}^{(d)}-1}\left(\Xi_{3\hat{p}^{(d)}-1}^{(d)}\right)$ of degree $3\hat{p}^{(d)}-1$ whose knot vector $\Xi_{3\hat{p}^{(d)}-1}^{(d)}$ follows \cref{eq:ReducedTripleProductKnotVector}, now applied to $\hat{\Xi}^{(d)}$, i.e.
\begin{equation}
        \label{eq:RefinedProjectionSpaceNumerics}
        \mathbb{S}_{\rho}^{(d)} = \mathbb{S}_{3\hat{p}^{(d)}-1}\left(\hat{\Xi}_{3\hat{p}^{(d)}-1}^{(d)}\right), \qquad \mathbb{S}_{\rho} = \mathbb{S}_{\rho}^{(1)} \otimes \mathbb{S}_{\rho}^{(2)} \otimes \mathbb{S}_{\rho}^{(3)}, 
\end{equation}
Since $\hat{\Xi}^{(d)}$ evolves from the geometry knot vector $\Xi^{(d)}$ by $p$- and $h$-refinement, it retains the regularity of $\Xi^{(d)}$ at the original knot values, and $\mathbb{S}_{3\hat{p}^{(d)}-1}\left(\Xi_{3\hat{p}^{(d)}-1}^{(d)}\right)$ therefore has, at these values, the same regularity as the univariate geometry space. At the same time, additional knots introduced by the $h$-refinement of the solution space are also contained in $\Xi_{3\hat{p}^{(d)}-1}^{(d)}$. The resulting projection space is therefore typically richer than the default space and is refined together with the solution discretization. In this sense, the construction can be regarded as an isoparametric choice for the projection space. \\

\begin{figure}[tp]
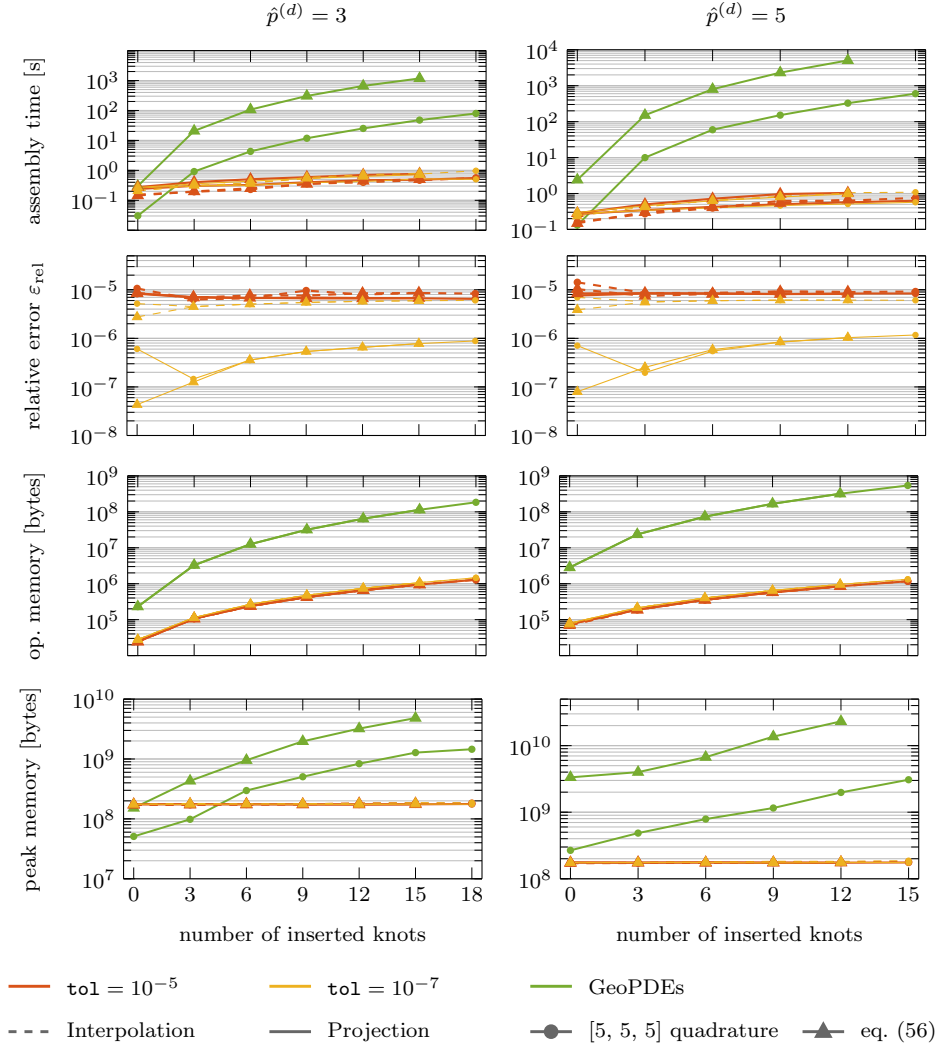

  \centering
  \setlength{\tabcolsep}{2pt}
  \begin{tabular}{cc}
    \input{flag_default_r1c1.tex} & \input{flag_default_r1c2.tex}\\
    \input{flag_default_r2c1.tex} & \input{flag_default_r2c2.tex}\\
    \input{flag_default_r3c1.tex} & \input{flag_default_r3c2.tex}\\
    \input{flag_default_r4c1.tex} & \input{flag_default_r4c2.tex}\\
  \end{tabular}\\[4pt]
  \input{flag_default_legend.tex}
  \caption{Results of the stiffness matrix/tensor assembly for the \emph{thick} flag geometry \cref{fig:1_a} using the default projection space \cref{eq:DefaultProjectionSpace} for $\hat{p}^{(d)} = 3$ and $\hat{p}^{(d)} = 5$ (columns).}
  \label{fig:flag_default_stiff}
\end{figure}

\begin{figure}[tp]
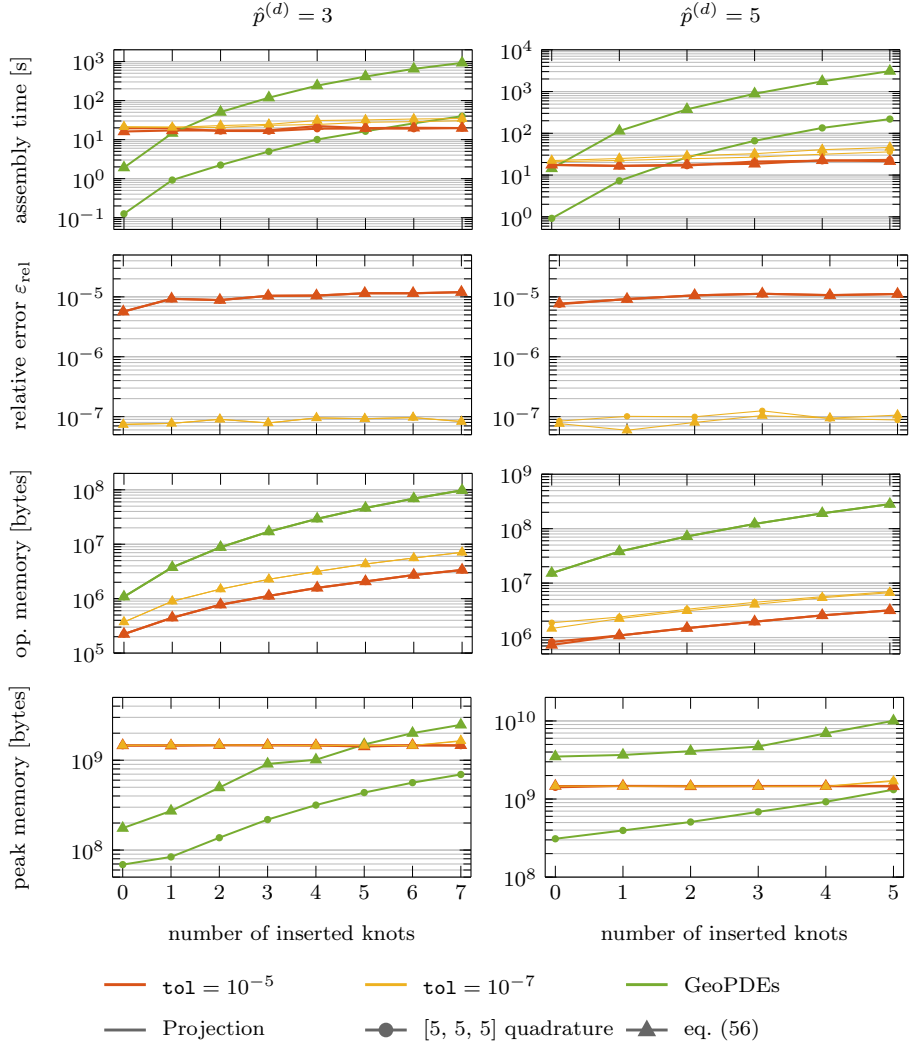

  \centering
  \setlength{\tabcolsep}{2pt}
  \begin{tabular}{cc}
    \input{beam_default_r1c1.tex} & \input{beam_default_r1c2.tex}\\
    \input{beam_default_r2c1.tex} & \input{beam_default_r2c2.tex}\\
    \input{beam_default_r3c1.tex} & \input{beam_default_r3c2.tex}\\
    \input{beam_default_r4c1.tex} & \input{beam_default_r4c2.tex}\\
  \end{tabular}\\[4pt]
  \input{beam_default_legend.tex}
  \caption{Results of the stiffness matrix/tensor assembly for the twisted beam geometry \cref{fig:1_c} using the default projection space \cref{eq:DefaultProjectionSpace} for $\hat{p}^{(d)} = 3$ and $\hat{p}^{(d)} = 5$ (columns).}
  \label{fig:beam_default_stiff}
\end{figure}

The results of the stiffness tensor assembly are summarized in \cref{fig:flag_default_stiff,fig:flag_refined_stiff} for the \emph{thick flag} geometry, in \cref{fig:rotor_default_stiff,fig:rotor_refined_stiff} for the rotor blade, in \cref{fig:singularity_default_stiff,fig:singularity_refined_stiff} for the linear almost singular geometry and in \cref{fig:beam_default_stiff,fig:beam_refined_stiff} for the twisted beam.\\

\begin{figure}[tp]
  \centering
  \setlength{\tabcolsep}{2pt}
  \begin{tabular}{cc}
    \input{flag_refined_r1c1.tex} & \input{flag_refined_r1c2.tex}\\
    \input{flag_refined_r2c1.tex} & \input{flag_refined_r2c2.tex}\\
    \input{flag_refined_r3c1.tex} & \input{flag_refined_r3c2.tex}\\
    \input{flag_refined_r4c1.tex} & \input{flag_refined_r4c2.tex}\\
  \end{tabular}\\[4pt]
  \input{flag_refined_legend.tex}
  \caption{Results of the stiffness matrix/tensor assembly for the \emph{thick} flag geometry \cref{fig:1_a} using the projection space \cref{eq:RefinedProjectionSpaceNumerics} derived from the solution space for $\hat{p}^{(d)} = 3$ and $\hat{p}^{(d)} = 5$ (columns).}
  \label{fig:flag_refined_stiff}
\end{figure}

The results of the stiffness tensor assembly for the \emph{thick flag} geometry are shown in \cref{fig:flag_default_stiff,fig:flag_refined_stiff}. Using the fixed projection space of \cref{eq:DefaultProjectionSpace}, the approximation error stagnates and, in particular for the tighter tolerance $\varepsilon=10^{-7}$, does not approach the prescribed threshold under further refinement. This indicates that the fixed auxiliary space is not sufficiently rich to represent the reciprocal determinant with the required accuracy. If the projection space is refined together with the discretization, the error decreases accordingly and approaches the prescribed tolerances. This improved accuracy comes at the cost of increasing assembly times and, in particular, a larger peak memory consumption due to the growing projection and interpolation systems. Nevertheless, the proposed projection-based assembly remains considerably faster than the interpolation-based approach while providing comparable or better accuracy. \\

\begin{figure}[tp]
  \centering
  \setlength{\tabcolsep}{2pt}
  \begin{tabular}{cc}
    \input{rotor_default_r1c1.tex} & \input{rotor_default_r1c2.tex}\\
    \input{rotor_default_r2c1.tex} & \input{rotor_default_r2c2.tex}\\
    \input{rotor_default_r3c1.tex} & \input{rotor_default_r3c2.tex}\\
    \input{rotor_default_r4c1.tex} & \input{rotor_default_r4c2.tex}\\
  \end{tabular}\\[4pt]
  \input{rotor_default_legend.tex}
  \caption{Results of the stiffness matrix/tensor assembly for the rotor geometry \cref{fig:1_b} using the default projection space \cref{eq:DefaultProjectionSpace} for $\hat{p}^{(d)} = 3$ and $\hat{p}^{(d)} = 5$ (columns).}
  \label{fig:rotor_default_stiff}
\end{figure}

\begin{figure}[tp]
  \centering
  \setlength{\tabcolsep}{2pt}
  \begin{tabular}{cc}
    \input{rotor_refined_r1c1.tex} & \input{rotor_refined_r1c2.tex}\\
    \input{rotor_refined_r2c1.tex} & \input{rotor_refined_r2c2.tex}\\
    \input{rotor_refined_r3c1.tex} & \input{rotor_refined_r3c2.tex}\\
    \input{rotor_refined_r4c1.tex} & \input{rotor_refined_r4c2.tex}\\
  \end{tabular}\\[4pt]
  \input{rotor_refined_legend.tex}
  \caption{Results of the stiffness matrix/tensor assembly for the rotor geometry \cref{fig:1_b} using the projection space \cref{eq:RefinedProjectionSpaceNumerics} derived from the solution space for $\hat{p}^{(d)} = 3$ and $\hat{p}^{(d)} = 5$ (columns).}
  \label{fig:rotor_refined_stiff}
\end{figure}

The rotor blade represents a more challenging geometry, which is reflected in the larger approximation errors and the increased computational effort. In particular, the fixed default projection space \cref{eq:DefaultProjectionSpace} is again not sufficiently rich to reach the tighter prescribed tolerances, whereas refining the projection space together with the solution space considerably improves the accuracy. For the refined projection space \cref{eq:RefinedProjectionSpaceNumerics}, the results also show that sufficiently accurate quadrature is required for a fair comparison. Since the assembly involves products of splines from increasingly rich spaces, the resulting polynomial degrees increase and a fixed low-order quadrature introduces an additional integration error that may dominate the actual approximation error. Using quadrature that is exact for the arising spline products removes this effect and reveals the expected accuracy of the projection. Although the refined projection space and the corresponding quadrature increase the time and memory requirements, the proposed projection-based method still outperforms the interpolation-based assembly for this more challenging geometry. \\

\begin{figure}[tp]
  \centering
  \setlength{\tabcolsep}{2pt}
  \begin{tabular}{cc}
    \input{singularity_default_r1c1.tex} & \input{singularity_default_r1c2.tex}\\
    \input{singularity_default_r2c1.tex} & \input{singularity_default_r2c2.tex}\\
    \input{singularity_default_r3c1.tex} & \input{singularity_default_r3c2.tex}\\
    \input{singularity_default_r4c1.tex} & \input{singularity_default_r4c2.tex}\\
  \end{tabular}\\[4pt]
  \input{singularity_default_legend.tex}
  \caption{Results of the stiffness matrix/tensor assembly for the almost singular geometry \cref{fig:1_d} using the default projection space \cref{eq:DefaultProjectionSpace} for $\hat{p}^{(d)} = 3$ and $\hat{p}^{(d)} = 5$ (columns).}
  \label{fig:singularity_default_stiff}
\end{figure}

\begin{figure}[tp]
  \centering
  \setlength{\tabcolsep}{2pt}
  \begin{tabular}{cc}
    \input{singularity_refined_r1c1.tex} & \input{singularity_refined_r1c2.tex}\\
    \input{singularity_refined_r2c1.tex} & \input{singularity_refined_r2c2.tex}\\
    \input{singularity_refined_r3c1.tex} & \input{singularity_refined_r3c2.tex}\\
    \input{singularity_refined_r4c1.tex} & \input{singularity_refined_r4c2.tex}\\
  \end{tabular}\\[4pt]
  \input{singularity_refined_legend.tex}
  \caption{Results of the stiffness matrix/tensor assembly for the almost singular geometry \cref{fig:1_d} using the projection space \cref{eq:RefinedProjectionSpaceNumerics} derived from the solution space for $\hat{p}^{(d)} = 3$ and $\hat{p}^{(d)} = 5$ (columns).}
  \label{fig:singularity_refined_stiff}
\end{figure}

For the stiffness tensor we additionally consider a geometry whose Jacobian determinant becomes small in part of the parameter domain. It is a tensor-product B-spline volume of degree $p=(1,1,1)$ with $\mathbf{n}=(2,2,2)$ on $[0,1]^3$, so that the univariate bases are the linear B-splines on the open knot vector $\left\{0,0,1,1\right\}$. For a deformation parameter $c>-1$, the control-point tensor $\mathbf{C}$ in \cref{eq:GeometryMapTensor} collects the eight vertices of the unit cube, with the single corner at $\hat{x}=(1,1,1)$ displaced along the first coordinate to $(1+c,1,1)$. This yields the geometry map
\begin{equation*}
        G\left( \hat{x}^{(1)}, \hat{x}^{(2)}, \hat{x}^{(3)} \right) = \left[ \left(1+c\,\hat{x}^{(2)}\hat{x}^{(3)}\right)\hat{x}^{(1)},\ \hat{x}^{(2)},\ \hat{x}^{(3)} \right]^\top.
\end{equation*}
Its Jacobian $\nabla G$ is upper triangular, so that $\det\left(\nabla G\left(\hat{x}\right)\right) = 1 + c\,\hat{x}^{(2)}\hat{x}^{(3)}$. Since $\hat{x}^{(2)}\hat{x}^{(3)}\in[0,1]$, the orientation-preserving condition \cref{eq:OrientationPreserving} holds for $c>-1$, and for $c\in(-1,0]$ the determinant ranges between $1+c$ and $1$. In particular, for $c=-1+10^{-5}$ we have $10^{-5}\le\det\left(\nabla G\left(\hat{x}\right)\right)\le 1$, so that the geometry remains regular but becomes nearly singular towards the parametric edge $\hat{x}^{(2)}=\hat{x}^{(3)}=1$, where the displaced corner approaches the vertex $(0,1,1)$. The geometry is therefore particularly interesting, since when $\omega\left(\hat{x}\right)=\det\left(\nabla G\left(\hat{x}\right)\right)$ becomes small, the entries $q_{kl}\left(\hat{x}\right)=N_{kl}\left(\hat{x}\right)/\omega\left(\hat{x}\right)$ of $Q$ can become very large.

For this experiment, we compare the stiffness tensors only with respect to the interior degrees of freedom. The entries associated with basis functions close to the nearly singular region become very large and are particularly difficult to approximate, so that they would otherwise dominate the relative error. The interpolation-based approach is especially affected by this behavior, since interpolation forces the spline approximation to attain the large values of $q_{kl}$ close to the nearly singular region, thereby deteriorating the approximation away from these points. The fixed default projection space is also not sufficiently rich to approximate the reciprocal determinant accurately. In contrast, using the refined projection space of \cref{eq:RefinedProjectionSpaceNumerics} considerably improves the approximation, and the error approaches the prescribed tolerance, with exact quadrature yielding the clearest convergence behavior. At the same time, the assembly time and memory requirements of the proposed method remain nearly constant under refinement and are substantially smaller than those of the interpolation-based approach. \\

\begin{figure}[tp]
  \centering
  \setlength{\tabcolsep}{2pt}
  \begin{tabular}{cc}
    \input{beam_refined_r1c1.tex} & \input{beam_refined_r1c2.tex}\\
    \input{beam_refined_r2c1.tex} & \input{beam_refined_r2c2.tex}\\
    \input{beam_refined_r3c1.tex} & \input{beam_refined_r3c2.tex}\\
    \input{beam_refined_r4c1.tex} & \input{beam_refined_r4c2.tex}\\
  \end{tabular}\\[4pt]
  \input{beam_refined_legend.tex}
  \caption{Results of the stiffness matrix/tensor assembly for the twisted beam geometry \cref{fig:1_c} using the projection space \cref{eq:RefinedProjectionSpaceNumerics} derived from the solution space for $\hat{p}^{(d)} = 3$ and $\hat{p}^{(d)} = 5$ (columns).}
  \label{fig:beam_refined_stiff}
\end{figure}

The twisted beam is again the most challenging of the considered geometries. We omit the interpolation-based method, since its interpolation system is singular for this geometry, as discussed above. Interestingly, the results in \cref{fig:beam_default_stiff} show that the fixed default projection space \cref{eq:DefaultProjectionSpace} is already sufficiently rich. The relative error follows the prescribed tolerance under refinement, so that the refined projection space of \cref{eq:RefinedProjectionSpaceNumerics} does not provide a relevant advantage for this example. At the same time, the assembly time and peak memory of the proposed method increase only mildly with the refinement, whereas both quantities grow considerably faster for the classical \textsc{GeoPDEs} assembly. This behavior suggests that the advantage of the proposed method becomes even more pronounced for finer discretizations.

\section{Conclusion}
\label{sec:Conclusion}

We have presented a projection-based low-rank method for assembling the mass and the stiffness tensor of orientation-preserving tensor-product B-spline geometries.

For the mass tensor, the construction exploits that the weight function $\omega = \det\left(\nabla G\right)$ is itself a piecewise polynomial and therefore contained in a reduced spline product space, whose coefficient tensor is obtained from the coordinate control-point tensors $\mathbf{C}^{(a)}$ by univariate transfer operators. With exact quadrature and without truncation, the representation \cref{eq:LowRankMassMatrixNew} is an exact reformulation of the Galerkin mass tensor.

For the stiffness tensor, we split the entries $q_{kl}$ of the matrix-valued weight function into the polynomial numerators $N_{kl}$, which are treated in the same way, and the reciprocal determinant $\rho = 1/\omega$, which is approximated by an $L^2$-projection onto a tensor-product spline space $\mathbb{S}_{\rho}$. All computations are carried out in the TT format, so that the coefficient tensors, of order nine and twelve before grouping, are never formed explicitly and the multivariate integrals reduce to univariate integrals and contracted products.

The numerical experiments indicate that the proposed method is competitive with the interpolation-based approach of~\cite{BuengerDolgovStoll:2020}. For the mass tensor, it is almost consistently faster and more accurate, and for the stiffness tensor it additionally shows a slight advantage in memory. It should be emphasized that the proposed method is of particular advantage in two situations. It remains applicable for the twisted beam in \cref{fig:1_c}, where the reduced regularity of the geometry makes the interpolation system \cref{eq:InterpolationSystem} singular, and it is more accurate for the nearly singular geometry in \cref{fig:1_d}, since the projection does not force the approximation to attain the large values of $q_{kl}$ close to the nearly singular region, where $\det\left(\nabla G\right) = 0$.

We point out that the accuracy of the stiffness assembly is governed by the richness of $\mathbb{S}_{\rho}$. The fixed default space \cref{eq:DefaultProjectionSpace} is not sufficient for all considered geometries, whereas the projection space \cref{eq:RefinedProjectionSpaceNumerics}, which is refined together with the solution space, follows the prescribed tolerances.

We further note that the orientation-preserving assumption \cref{eq:OrientationPreserving} can be relaxed to a determinant of constant sign. If $\det\left(\nabla G\left(\hat{x}\right)\right) < 0$ for all $\hat{x} \in [0,1]^3$, then $\left|\det\left(\nabla G\right)\right| = -\det\left(\nabla G\right)$, and it suffices to negate the coefficient tensor $\tilde{\mathbf{C}}_{\Sigma}$ in \cref{eq:WeightFunctionFrobeniusProductCoefficientsTT}, which is a sign change in a single TT core. Both the projection system \cref{eq:ReciprocalProjectionLinearSystem} and the assembled operators are then obtained as before. What is essential is that the sign does not change inside the parametric domain.

As an outlook, the following questions remain open. A reasonable rounding heuristic for the assembly of both the mass and the stiffness tensor has to be found, so that a prescribed tolerance is attained at the lowest possible TT ranks. Moreover, the cost of the proposed stiffness assembly is dominated by the univariate operators associated with $\mathbb{S}_{\rho}$ and by the resulting projection system \cref{eq:ReciprocalProjectionLinearSystem,eq:TransformedLinearSystem}, which is why projection spaces that approximate $\rho$ with fewer basis functions, and preconditioners adapted to the weighted projection system, appear to be the most promising direction for further improvement. Further, an extension of the proposed method to multi-patch geometries, cf.~\cite{Riemer2025}, and to hierarchical discretizations, cf.~\cite{Riemer2026}, is of interest. Finally, the assembled TT operators could be used in low-rank solvers, cf.~\cite{montardini2023lowsp,Mika2026}.


%
\section*{Conflict of interest}
The authors declare that they have no conflict of interest.

The authors did not receive support from any organization for the submitted work.

The authors have no competing interests to declare that are relevant to the content of this article.

The computer codes and algorithms used during the current study are available at \url{https://github.com/TomRiem/Low_Rank_IGA_Projection.git}.



\bibliographystyle{spmpsci}
\bibliography{literature}

\end{document}